\documentclass[10pt,oneside,a4paper]{article} 						

\usepackage[fleqn]{amsmath} 										
\usepackage{amsthm,amsfonts,latexsym,amssymb,amscd} 				
\usepackage[mathscr]{eucal}   										
\usepackage[all,2cell]{xy} \UseAllTwocells 							
\usepackage{stmaryrd} 												
\usepackage{framed}
\usepackage{color}
\usepackage{txfonts}
\usepackage{stackrel}
\usepackage{scalerel} 	
\usepackage{cite}

\usepackage{tikz}
\usetikzlibrary{svg.path}

\definecolor{orcidlogocol}{HTML}{A6CE39}
\tikzset{
  orcidlogo/.pic={
    \fill[orcidlogocol] svg{M256,128c0,70.7-57.3,128-128,128C57.3,256,0,198.7,0,128C0,57.3,57.3,0,128,0C198.7,0,256,57.3,256,128z};
    \fill[white] svg{M86.3,186.2H70.9V79.1h15.4v48.4V186.2z}
                 svg{M108.9,79.1h41.6c39.6,0,57,28.3,57,53.6c0,27.5-21.5,53.6-56.8,53.6h-41.8V79.1z M124.3,172.4h24.5c34.9,0,42.9-26.5,42.9-39.7c0-21.5-13.7-39.7-43.7-39.7h-23.7V172.4z}
                 svg{M88.7,56.8c0,5.5-4.5,10.1-10.1,10.1c-5.6,0-10.1-4.6-10.1-10.1c0-5.6,4.5-10.1,10.1-10.1C84.2,46.7,88.7,51.3,88.7,56.8z};
  }
}

\newcommand\orcidicon[1]{\href{https://orcid.org/#1}{\mbox{\scalerel*{
\begin{tikzpicture}[yscale=-1,transform shape]
\pic{orcidlogo};
\end{tikzpicture}
}{|}}}}

\usepackage[draft=false,breaklinks,backref]{hyperref} 					

\newcommand{\E}{\mathcal{E}}

\newcommand{\Fg}{\mathfrak{F}}\newcommand{\Gg}{\mathfrak{G}}

\newcommand{\CC}{{\mathbb{C}}}

\newcommand{\FF}{{\mathbb{F}}}
 
\newcommand{\KK}{{\mathbb{K}}}

\newcommand{\NN}{{\mathbb{N}}}
\newcommand{\QQ}{{\mathbb{Q}}}
\newcommand{\RR}{{\mathbb{R}}} 

\newcommand{\ZZ}{{\mathbb{Z}}}

\newcommand{\As}{{\mathscr{A}}}\newcommand{\Bs}{{\mathscr{B}}}\newcommand{\Cs}{{\mathscr{C}}}
\newcommand{\Ds}{{\mathscr{D}}}

\newcommand{\Is}{{\mathscr{I}}}\newcommand{\Js}{{\mathscr{J}}} 
\newcommand{\Ms}{{\mathscr{M}}}\newcommand{\Ns}{{\mathscr{N}}}
\newcommand{\Ps}{{\mathscr{P}}}	
\newcommand{\Qs}{{\mathscr{Q}}}\newcommand{\Rs}{{\mathscr{R}}}
\newcommand{\Ts}{\mathscr{T}} 			
\newcommand{\Vs}{{\mathscr{V}}}

\newcommand{\Zs}{{\mathscr{Z}}}

\DeclareFontFamily{U}{rsfs}{\skewchar\font127 }
\DeclareFontShape{U}{rsfs}{m}{n}{%
   <5> <6> rsfs5
   <7> rsfs7
   <8> <9> <10> <10.95> <12> <14.4> <17.28> <20.74> <24.88> rsfs10
}{}
\DeclareSymbolFont{rsfs}{U}{rsfs}{m}{n} 
\DeclareSymbolFontAlphabet{\scr}{rsfs}

\newcommand{\Af}{\scr{A}}\newcommand{\Bf}{\scr{B}}
\newcommand{\If}{\scr{I}}
\newcommand{\Mf}{\scr{M}} 
\newcommand{\Sf}{\scr{S}}

\DeclareMathOperator{\ke}{Ker}
\DeclareMathOperator{\im}{Im}
 
\DeclareMathOperator{\id}{Id} 

\DeclareMathOperator{\Ob}{Ob}

\DeclareMathOperator{\Hom}{Hom}
 
\DeclareMathOperator{\dom}{Dom}

\DeclareMathOperator{\ev}{ev}

\DeclareMathOperator{\End}{End}

\DeclareMathOperator{\Diff}{Diff}

\renewcommand{\emph}{\textbf} 										
\newcommand{\ip}[2]{\langle #1\mid #2\rangle}						
\renewcommand{\iff}{\Leftrightarrow}								
\newcommand{\imp}{\Rightarrow}										
\newcommand{\st}{\ : \ }											

\newcommand{\hlink}[2]{\href{#1}{\texttt{#2}}} 						

\newcommand{\xqedhere}[2]{%
  \rlap{\hbox to#1{\hfil\llap{\ensuremath{#2}}}}}

\newcommand{\xqed}[1]{%
  \leavevmode\unskip\penalty9999 \hbox{}\nobreak\hfill
  \quad\hbox{\ensuremath{#1}}}

\theoremstyle{plain}
\newtheorem{theorem}{Theorem}[section]							

\newtheorem{proposition}[theorem]{Proposition}

\newtheorem{definition}[theorem]{Definition}

\theoremstyle{definition} 
\newtheorem{remark}[theorem]{Remark}

\numberwithin{equation}{section}  									
\title{\textbf{Duality for Multimodules}}

\author{\normalsize 
\orcidicon{0000-0002-1387-9283} 
Paolo Bertozzini $^a$  
\ 
\orcidicon{0000-0003-3128-6884}
Roberto Conti $^b$
\ 
\orcidicon{0000-0002-1451-5000}
Chatchai Puttirungroj $^c$   
\\  
\normalsize $^{a \ c}$ \textit{Department of Mathematics and Statistics, Faculty of Science and Technology,}
\\
\normalsize \textit{Thammasat University, Pathumthani 12121, Thailand}
\\
\normalsize e-mail: $^a$ \texttt{paolo.th@gmail.com} \quad $^c$ \texttt{cputtirungroj@gmail.com} 
\\ 
\normalsize $^b$ \textit{Dipartimento di Scienze di Base e Applicate per l'Ingegneria, 
} 
\\
\normalsize \textit{Sapienza Universit\`a di Roma, Via A. Scarpa 16, I-00161 Roma, Italy}
\\
\normalsize e-mail: \texttt{roberto.conti@sbai.uniroma1.it} 
}

\date{\normalsize{02 August 2021}}

\begin{document}

\maketitle

\begin{center}
\textit{dedicated to N.Bourbaki who initiated the study of multimodules
} 
\end{center} 

\begin{abstract} \noindent 
Suitable duals of multimodules are introduced and used to provide transposition contravariant right semi-adjunctions (and dualitites under reflexivity). 
Several additional notions on multimodules are discussed: generalized morphisms and involutions; tensor products and contractions; inner products; first-order differential operators. A construction of an (involutive) colored properad of multimodules is suggested. 

\medskip


\smallskip

\noindent
\emph{MCS-2020:} 
					16D10,			
					16D80, 			
					16D90,			
					16D99, 			
					18M85. 			
\end{abstract}

\tableofcontents


\newpage

\section{Introduction and Motivation} \label{sec: intro}

The study of categories of modules and bimodules over unital associative rings or algebras is one of the most developed subjects of modern algebra and its inception might be traced back to the work of R.Dedekind, E.Noether and B.L.van den Waerden, among many others. 

\medskip 

Multimodules over unital associative rings and algebras are quite a natural generalization of right/left-modules and bimodules that, as far as we know, have been first described in N.Bourbaki~\cite{Bou89}. Since, for most purposes, multimodules are equivalently seen as bimodules over tensor products of rings and algebras, it can be claimed that their investigation essentially reduces to the study of special classes of bimodules and not much attention has been paid to them (we have been able to locate only one specific reference on multimodules~\cite{Ke62} and some sporadic mentioning of them, for example in~\cite[section~0]{Ta87}). 

\medskip 

The ``substitution'' of multimodules with corresponding bimodules over tensor products turns out to be problematic whenever the category of morphisms is extended with the inclusion of maps that have different covariance properties with respect to the several actions involved. One could still substitute multimodules with bimodules over tensor products of rings, as long as such tensor products of rings are simultaneously equipped with different products (all distributive with respect to the same Abelian group structure), but this essentially amounts to define an ``hyper-algebra structure'' on the tensor product multimodule of the rings (see remark~\ref{rem: hyper}). 

\medskip 

The basic algebraic material here presented naturally arose as a byproduct in our study of non-commutative generalizations of contravariant calculus.~\footnote{
P.Bertozzini, R.Conti, C.Puttirungroj, Non-commutative Contravariant Differential Calculus (in preparation). \label{foo: ncdc}
} 
Since quite surprisingly we have not been able to locate any relevant source dealing with this topic, we thought that the subject deserves an adequate separate treatment. 
Specifically (anticipating arguments and motivations pertaining to the aforementioned work) in non-commutative (algebraic) geometry, it is a common thread to look for generalizations of the usual notion of ``differential operator'' to the case of maps between bimodules over a non-commutative algebra $\As$ and it often happens (for example whenever one is considering ``double derivations'' on $\As$) that the spaces of such ``non-commutative differential operators'' are naturally equipped with a multimodule structure over the original non-commutative algebra $\As$. 
Although of tangential interest for this work, a general definition of first-order differential operators between multimodules, covering in particular all such cases, will be included in appendix~\ref{sec: 1st-ord-multi}. 
Further developments in these directions, including investigations of non-commutative vector fields and non-commutative connections on (multi-)modules, will have to wait subsequent works (see the paper in footnote~\ref{foo: ncdc} and references therein for more details). 

\bigskip 

In short, the specific goals of the present work are to:  
\begin{itemize} \it
\item[$\rightsquigarrow$] 
define multimodules based over an arbitrary $\Zs$-central bimodule~\footnote{
Where $\Zs$ is a commutative unital associative ring/algebra} (more generally over a \hbox{$\Zs$-central} unital associative ring $\Rs_\Zs$) instead of just an Abelian group:  
\end{itemize}
this allows to discuss mutually commuting (right/left) actions that are compatible with a certain fixed $\Zs$-linear structure, but that can still have alternative $\Rs$-linearity properties;  
\begin{itemize} \it
\item[$\rightsquigarrow$] 
introduce a notion of involution for multimodules that allows for different covariance/contravariance:
\end{itemize}
since involutions for us are just involutive morphisms, this requires an appropriate definition of category of multimodules, where morphisms (necessarily $\Zs$-linear) can have different covariance properties (and even different conjugate-$\Rs_\Zs$-linearity properties) with respect to the different actions involved; 
\begin{itemize}\it
\item[$\rightsquigarrow$] 
provide a systematic treatment of the several ($\Zs$-linear) duals of multimodules, their associated categorical semi-adjunctions and (under saturation conditions for evaluations) establish transposition dualities: 
\end{itemize}
it is already known that in the case of bimodules one needs to separately consider right and left duals in place of the usual notion of dual vector space; in the case of multimodules, the situation is a bit more involved and one can construct different (conjugate)-duals for any choice of subfamilies of left/right actions (and corresponding conjugations of $\Rs_\Zs$); each dual is defined in this work via a universal factorization property and its elements are concretely realized as $\Zs$-multilinear functions that have selective \hbox{$\Rs_\Zs$-(conjugate)}-linearity properties with respect to the specified actions; 
\begin{itemize}
\item[$\rightsquigarrow$] \it 
introduce universal traces, and more generally contractions, on multimodules; 
\end{itemize}
traces of linear operators and contractions of tensors are quite standard operations performed in multilinear algebra; we reframe such notions in the more general context of multimodules, providing again a definition via universal factorization properties; 
\begin{itemize}
\item[$\rightsquigarrow$] \it 
discuss, for multimodules over involutive rings/algebras, suitable notions of ``inner products'' and (under conditions of non-degeneracy/fullness) establish Riesz isomorphisms: 
\end{itemize}
inner products on multimodules also come in several types, each corresponding to a different dual, and are here realized as certain \textit{balanced multi-sesquilinear maps}; involutive algebras are necessary in order to give a meaning to Hermitianity conditions on inner products; every inner-product induces a \textit{canonical Riesz morphisms} of a multimodule  into a corresponding dual; non-degeneracy and fullness are required to obtain an isomorphism. 
\begin{itemize}\it
\item[$\rightsquigarrow$] 
describe first order differential operators between multimodules: 
\end{itemize}  
the first order condition in non-commutative geometry~\cite[sections $4.\gamma$ and $4.\delta$]{Co94}, usually formulated in the case of operators between bimodules, is here expanded to cover the general setting of multimodules; 
\begin{itemize}
\item[$\rightsquigarrow$] \it
make the first steps toward a study of involutive colored properads using multimodules as a template: 
\end{itemize}
the material here included is mostly intended to provide a usable language for quite practical situations (some of which have been actually originating from work in categorical non-commutative geometry) where multimodules and their duals might be used and manipulated. As a consequence, we have not been looking for maximal generality in the statements and we kept a rather low sophistication level in the discussion of all the category-theoretical aspects of the subject; 
a more detailed study of these topics is under way, but we can already anticipate that it will fall within the scope of certain variants of \textit{involutive colored properads} and \textit{involutive polycategories}. 
As stated above, we plan to address more properly these points in subsequent works. 

\bigskip 

Here below is a more detailed description of the content of the paper. 

\medskip 

In section~\ref{sec: generalities} we modify the usual setting of bimodules over unital associative rings considering, in place of the initial ring $\ZZ$ a commutative unital associative ring $\Zs$ and instead of rings acting on Abelian groups (\hbox{$\ZZ$-bimodules}), $\Zs$-central algebras $\As$ acting in a $\Zs$-bilinear way on $\Zs$-central bimodules $\Ms$. 
Morphisms are in this case pairs of $\Zs$-linear maps (in place of additive maps) that induce a unital $\Zs$-linear covariant or contravariant grade-preserving homomorphism on the associated $\NN$-graded algebras $\Ms^\urcorner:=\As\oplus\Ms\oplus\{0\}\cdots$ of the bimodules. 
This kind of environment can immediately describe, as a special case, categories of $\KK$-linear covariant or contravariant morphisms of unital bimodules over $\KK$-algebras, for a certain field $\KK$ in place of $\Zs$.  

\medskip 

The existence of many situations requiring the usage of non-trivial (involutive auto)morphisms for the base field $\KK$ and the consequent need to deal simultaneously with maps that are not $\KK$-linear, imposes a further refinement of the structure: the common base commutative associative unital ring $\Zs$ is replaced by a \hbox{$\Zs$-central} unital associative ring $\Rs_\Zs$. The family of unital covariant or contravariant $\Zs$-linear homomorphisms $\phi$ of $\Zs\oplus\Rs$ identifies the possible alternative notions of $\phi$-linearity with respect to the base ring $\Rs$. 
The paradigmatic situation with $\Rs:=\CC$ and $\Zs:=\RR$ imposes only $\RR$-linearity on morphisms that are further classified as $\CC$-linear and $\CC$-conjugate-linear depending on the choice of the $\CC$-automorphism $\phi$; but the formalism can be used in the case of algebras over arbitrary extensions of fields (or more generally extensions of rings). 

\medskip 

Section~\ref{sec: multimodules} puts forward our definition of multimodules over families of unital associative algebras over $\Rs_\Zs$. We stress that taking $\Rs=\Zs=\ZZ$, we just reproduce the usual definition of multimodules in~\cite{Bou89} and taking $\Zs\to\Rs$ an extension of fields we obtain multimodules as $\Zs$-vector spaces equipped with $\Zs$-bilinear actions of $\Rs$-algebras, where morphisms can be $\phi$-linear for any $\Zs$-linear automorphism $\phi$ of $\Rs$. 
The unavoidability of multimodules (in every context dealing with bimodules) is witnessed by the construction of $\Zs$-central multimodules of $\Zs$-linear maps, and $\Zs$-tensor products, between $\Zs$-central bimodules. 

\medskip 

In section~\ref{sec: involutions-m} we specialize to the treatment of involutive endomorphisms of $\Zs$-central multimodules over \hbox{$\Zs$-central} $\Rs_\Zs$-algebras and we examine how involutions on bimodules (and multimodules) propagate to involutions for spaces of $\Zs$-linear morphisms and tensor products of multimodules. 

\medskip 

The main result of the paper is contained in section~\ref{sec: multi-dual-paring} where we introduce definitions of duals of multimodules via universal factorization properties and we prove that transpositions functors in the category of  multimodules give rise to contravariant right semi-adjunctions (theorem~\ref{th: sadj}) that, for multimodules satisfying reflexivity, produce dualities. 
In general there exist different \textit{conjugate-duals} for a $\Zs$-central multimodule ${}_{(\As_\alpha)_A}\Ms_{(\Bs_\beta)_B}$ over \hbox{$\Rs_\Zs$-algebras}, each one of them ${}^{(\gamma_i)_I}\Ms^{(\gamma_j)_J}$ specified by certain families $(\gamma_i)_{i\in I}$ and $(\gamma_j)_{j\in J}$ of $\Rs$-conjugations, with arbitrary sets of indexes $I\subset A$ and $J\subset B$.

\medskip 

In the first part of section~\ref{sec: traces-ip} we define universal contractions/traces on multimodules via universal factorization properties and we construct them quotienting the original multimodules with respect to certain \textit{commutator sub-multimodules}. 
The remaining part of section~\ref{sec: traces-ip} discusses tentative generalizations, to the setting of multimodules over involutive algebras, of the familiar notion of inner product for vector spaces or modules and for each such inner product defines its Riesz ``natural transformation''.\footnote{
Due to the different covariance of the functors involved, a categorical discussion of the ``naturality'' of Riesz morphisms would require the usage of hybrid 2-categories~\cite{BePu14} of multimodules. 
} 
Under conditions of non-degeneracy and fullness of the inner products, we also provide a multimodule version of Riesz isomorphism theorem.  
The inner products here introduced are not necessarily positive: a positivity requirement can be added (at such an abstract level) imposing the existence of positive cones on the algebras. 

\medskip 

The final outlook section~\ref{sec: outlook} briefly expands on the already mentioned planned utilization of the categories of multimodules, here developed, as a paradigmatic example in the study of the abstract notion of ``involutive colored properad'' and their associated involutive ``convolution hyper-algebroids'' following the lines that some of us have discussed in previous papers~\cite{BCLS20}. 

\medskip 

In appendix~\ref{sec: functorial-pairing}, we briefly recall the notion of (contravariant) semi-adjunction~\cite{Med74}, a special case of regular full functorial pairings later defined in~\cite{Wis13}, that will be needed to describe the dualities for contravariant trasposition functors in categories of multimodules.
Special attention has been devoted to the explicit characterization of semi-adjunctions for contravariant functors. 

\medskip

As already mentioned, the present paper was motivated by an ongoing effort towards the study of non-com\-mu\-ta\-tive vectors fields and contravariant non-commutative differential calculus (see footnote~\ref{foo: ncdc}); in appendix~\ref{sec: 1st-ord-multi} we present the generalization, to the case of multimodules, of a definition of first-order differential operator on bimodules over non-commutative algebras, that has been useful in that context. Further extensions in the direction of differential analysis on multimodules (starting with a theory of connections) are briefly mentioned in the outlook section and will be dealt with elsewhere.  

\section{Generalities} \label{sec: generalities}

We start specifying basic settings and definitions; for more details on background material that is not explicitly mentioned, we refer to the texts~\cite{Al09} and~\cite{Bou89}. 

\medskip 

We assume $\Zs$ to be a commutative unital associative ring. 
All the rings $\Rs$ here considered will be unital associative (not necessarily commutative) and \emph{$\Zs$-central rings}: they are equipped with a unital homomorphism of rings $\iota_\Rs:\Zs\to Z(\Rs):=\{r\in\Rs \ | \ \forall x\in\Rs \st r\cdot x=x\cdot r \}$, where $Z(\Rs)$ denotes  the center of the ring $\Rs$ (itself a commutative unital associative ring). 

\medskip 

All the $\Rs$-bimodules $\Ms$ considered in this paper are assumed to be unital ($1_\Rs\cdot x=x$, for all $x\in\Ms$) and 
\emph{\hbox{$\Zs$-central} $\Rs$-bimodules}, meaning that there is a unital homomorphism of rings $\iota_\Ms:\Zs\to Z(\Ms)^0$, where we define $Z(\Ms)^0:=\{r\in Z(\Rs) \ | \ \forall x\in\Ms \st r\cdot x=x\cdot r \}$ as  
the center ring of the $\Rs$-bimodule $\Ms$ (that is itself a $\Zs$-central unital sub-ring of $\Rs$). Similarly, 
$Z(\Ms)^1:=\{x\in\Ms \ | \ \forall r\in\Rs \st r\cdot x=x\cdot r \}$ denotes the center module of $\Ms$ (itself a $\Zs$-central unital $\Rs$-bimodule). 

\medskip 

Here \emph{$\Zs$-central $\Rs$-algebras} are defined as $\Zs$-central $\Rs$-bimodules $\As:={}_\Rs\As_\Rs$ with a distributive multiplication $\circ$ such that: 
$(r\cdot x)\circ y=r\cdot(x\circ y)$, $(x\cdot r)\circ y=x\circ(r\cdot y)$ and $x\circ(y\cdot r)=(x\circ y)\cdot r$, for all $x,y\in\As$ and $r\in \Rs$. In this way, multiplication in a $\Zs$-central $\Rs$-algebra is necessarily $\Zs$-bilinear and every $\Zs$-central ring $\Rs$ becomes an (associative unital) \hbox{$\Zs$-central} algebra over itself.  
We will usually consider $\Zs$-central $\Rs$-algebras $\As$ that are unital and associative. We will consider $\Zs$-central $\As$-bimodules $\Ms$ that are unital and hence become canonically $\Zs$-central $\Rs$-bimodules with action $r\cdot x:=(r\cdot 1_\As)\cdot_\Ms x$, for $r\in\Rs$ and $x\in\Ms$.

\medskip 

Since $\ZZ$ is initial in the category of unital associative rings, $\ZZ\cdot 1_\Rs\subset \Zs\subset Z(\Rs)$ and the \emph{characteristic of $\Rs$} is the minimum $n\in\NN$ such that $\ke(\iota)=n\cdot\ZZ$, where $\iota:\ZZ\to\Rs$ is the initial unital homomorphism $z\mapsto z\cdot 1_\Rs$. Whenever the characteristic is a prime number, $\Rs$ is actually an $\FF$-algebra over the finite field $\FF:=\ZZ/\ke(\iota)$. 

\medskip 

Particular attention should be given to the definition of morphisms for bimodules over $\Zs$-central $\Rs$-algebras. 
\begin{definition}
Let $\Zs$ be a commutative unital associative ring and $\Rs_\Zs$ a unital associative $\Zs$-central ring.  

A map $\Ms\xrightarrow{\Phi}\Ns$ between two $\Zs$-central unital bimodules $\Ms:=\Ms_\Zs$, $\Ns:=\Ns_\Zs$,  
is said to be \emph{$\Zs$-linear} if:
\begin{equation*} 
\Phi(x+y)=\Phi(x)+\Phi(y), \quad \quad \Phi(\iota_\Ms(z)\cdot x)=\iota_\Ns(z)\cdot\Phi(x),\quad  \forall x,y\in\Ms, \ z\in\Zs.
\end{equation*}

A map $\As\xrightarrow{\phi}\Bs$ between $\Zs$-central unital associative rings $\As:=\As_\Zs$, $\Bs:=\Bs_\Zs$ is  
\begin{itemize}
\item
\emph{covariant} if: $\phi(x \circ_\As y)=\phi(x)\circ_\Bs\phi(y)$, for all $x,y\in\As$, 
\item 
\emph{contravariant} if: $\phi(x \circ_\As y)=\phi(y)\circ_\Bs\phi(x)$, for all $x,y\in\As$, 
\item
\emph{unital} if: $\phi(1_\As)=1_\Bs$, 
\item
\emph{homomorphism} if: it is $\Zs$-linear covariant and unital, 
\item
\emph{anti-homomorphism} if: it is $\Zs$-linear contravariant and unital. 
\end{itemize}

A $\Zs$-linear map $\Ms:={}_\As\Ms_{\As}\xrightarrow{\Phi}{}_\Bs\Ns_{\Bs}$ between $\Zs$-central unital bimodules over $\Zs$-central unital associative rings $\As:=\As_\Zs,\Bs:=\Bs_\Zs$ is said to be \emph{$\phi$-linear}, for a certain $\Zs$-linear unital homomorphism (anti-homomorphism) $\As\xrightarrow{\phi}\Bs$ if: 
\begin{align*}
&\Phi(a_1\cdot x\cdot a_2)=\phi(a_1)\cdot \Phi(x)\cdot \phi(a_2), \quad \forall x\in\Ms, \ a_1,a_2\in\As, 
&&\text{in the $\phi$-covariant case}, 
\\
&\Phi(a_1\cdot x\cdot a_2)=\phi(a_2)\cdot \Phi(x)\cdot \phi(a_1), \quad \forall x\in\Ms, \ a_1,a_2\in\As,
&&\text{in the $\phi$-contravariant case}.
\end{align*}

A $\phi$-linear covariant (contravariant) \emph{morphism} of $\Zs$-central unital bimodules, over $\Zs$-central unital associative rings, consists of a pair $(\phi,\Phi)$ as above. In the case of $\Zs$-central unital associative algebras over 
$\Zs$-central unital associative rings, the morphism $\As\xrightarrow{\Phi}\Bs$ must be unital and covariant (contravariant).  

\medskip 

For $\Zs$-central bimodules ${}_\As\Ms_\As, {_\Bs}\Ns_\Bs$ over $\Zs$-central unital associative algebras ${}_\Rs\As_\Rs, {}_\Rs\Bs_\Rs$ over a $\Zs$-central unital associative ring $\Rs_\Zs$, morphisms are still denoted by $(\phi,\Phi)$, where $\Phi:=(\Phi^0,\Phi^1)$ is a pair of $\phi$-linear unital morphisms $\Phi^0:\As\to\Bs$ of algebras and $\Phi^1:\Ms\to\Ns$ of bimodules, such that $\Phi^1$ is $\Phi^0$-linear: 
\begin{equation*}
\xymatrix{
\Zs \ar@{=}[d] \ar[r]^{\iota_\Rs} & \Rs \ar[d]^\phi \ar[r]^{\iota_\As} & \As\ar[d]^{\Phi^0} & \Ms \ar[d]^{\Phi^1}
&\Phi^1(a_1\cdot x\cdot a_2)=\Phi^0(a_1)\cdot \Phi^1(x)\cdot \Phi^0(a_2)
\quad 
\text{covariant $\phi$-linear case\phantom{ntra}} 
\\
\Zs \ar[r]^{\iota_\Rs} & \Rs \ar[r]^{\iota_\Bs} & \Bs  & \Ns 
&  
\Phi^1(a_1\cdot x\cdot a_2)=\Phi^0(a_2)\cdot \Phi^1(x)\cdot \Phi^0(a_1)
\quad 
\text{contravariant $\phi$-linear case}. 
}
\end{equation*}
\end{definition}

\begin{remark}
We have a category of $\Zs$-linear maps between $\Zs$-central unital bimodules over $\Zs$-central unital associative $\Rs_\Zs$-algebras. Such category is not $\ZZ_2$-graded with respect to covariance / contravariance, since the same morphism $\Phi$ can be $\phi$-covariant or $\phi$-contravariant depending on the choice of $\phi$. 

\medskip 

A better alternative consists, as we did, in defining morphisms as triples ${}_\As\Ms_\As\xrightarrow{(\phi,\Phi)}{}_\Bs\Ns_\Bs$ of $\Zs$-linear maps $\phi:\Rs_\Zs\to\Rs_\Zs$, $\Phi^0:\As_\Rs\to\Bs_\Rs$ and $\Phi^1:\Ms_\As\to\Ns_\Bs$ with $(\phi,\Phi^0)$ and $(\Phi^0,\Phi^1)$ both $\phi$-linear morphisms. 
In this case the category is $\ZZ_2$-graded (by the covariance of the triple) furthermore it is isomorphic to the \hbox{$\ZZ_2$-graded} category of degree zero unital $\Zs$-linear (covariant or contravariant) morphisms $(\phi,\Phi^0,\Phi^1)$ between graded unital associative \hbox{$\Zs$-central} algebras of the form $\Ms^\urcorner:=\Rs\oplus\As\oplus \Ms\oplus\{0\}\cdots$. 
\xqed{\lrcorner}
\end{remark}

\begin{definition}
A \emph{covariant} (respectively \emph{contravariant}) \emph{involution} on a $\Zs$-central unital associative ring $\Rs_\Zs$ is a $\Zs$-linear covariant (respectively contravariant) map $\Rs\xrightarrow{\star}\Rs$ that is involutive $(x^\star)^\star=x$, for all $x\in\Rs$. 

\medskip 

Whenever dealing with $\Zs$-central algebras $\As_\Rs$ over a $\Zs$-central unital associative ring $\Rs_\Zs$, we use the term \emph{$\Rs_\Zs$-conjugation} to denote an involution of the $\Zs$-central unital associative ring $\Rs_\Zs$. 

\medskip 

A \emph{covariant (contravariant) involution} $\star$ on ${}_\Rs\As_\Rs$ is said to be \emph{$\gamma$-conjugate-linear} if it is $\gamma$-linear for a certain covariant (contravariant) $\Rs_\Zs$-conjugation $\gamma$, specifically:  
$(r_1\cdot x\cdot r_2)^\star=\gamma(r_1)\cdot x^\star\cdot \gamma(r_2)$ in the $\gamma$-covariant case;  
$(r_1\cdot x\cdot r_2)^\star=\gamma(r_2)\cdot x^\star\cdot \gamma(r_1)$) in the $\gamma$-contravariant case, for all $r_1,r_2\in\Rs$ and $x\in\As$. 
\end{definition}

\begin{remark}
For $\Zs$-central algebras $\As_\Rs$ over non-commutative rings $\Rs$, covariant (contravariant) involutions can be $\gamma$-conjugate-linear only with respect to a covariant (contravariant) conjugation $\gamma$. 
Whenever $\Rs$ is commutative, there is no difference between covariant and contravariant conjugations and hence, for an arbitrary conjugation $\gamma$, we can have covariant or contravariant involutions on $\As_\Rs$ that are $\gamma$-conjugate-linear. 

\medskip 

Notice that for involutive $\Zs$-central $\Rs$-algebras $\As$, we necessarily have $\eta_\As(\Zs)\subset Z(\As_\Rs)^0\cap\{x\in\As \ | \ x^\star=x \}$. It is of course possible, for a certain $\Zs$-central ring $\Rs$ to have involutions $\gamma$ that do not necessarily leave $\eta_\Rs(\Zs)$ invariant or that do not necessarily fix all the elements of $\eta_\Rs(\Zs)$; in this case one can further ``restrict'' the commutative algebra $\Zs$ in order to make $\gamma$ a conjugation: given a certain family $\Gamma$ of additive (covariant or contravariant) involutions of $\Rs$, we see that 
$\Zs^\Gamma:=\eta^{-1}_\Rs\left(\bigcap_{\gamma\in\Gamma} \{x\in Z(\Rs) \ | \ \gamma(r)=r \}\right)$ is a unital sub-algebra of $\Zs$ making all the $\gamma\in\Gamma$ conjugations of $\Rs$ as $\Zs^\Gamma$-central ring. 
\xqed{\lrcorner}
\end{remark}

\medskip 

There are universal ways to reformulate $\gamma$-conjugate-linear unital morphisms of $\Zs$-central $\Rs$-algebras (and also of $\Zs$-central $\Rs$-bimodules) as covariant $\Rs$-linear unital morphisms. 
\begin{definition}\label{def: cj-dual}
Given a $\Zs$-central $\Rs$-algebra $\As_\Rs$ and a conjugation $\gamma$ in $\Rs$, a \emph{$\gamma$-conjugate of $\As_\Rs$} 
consists of a $\gamma$-conjugate-linear unital morphism of $\Zs$-central $\Rs$-algebras $\As\xrightarrow{\eta_\As}\As^\gamma$ that satisfies the universal factorization property: 
for any $\gamma$-conjugate-linear unital morphism of $\Zs$-central $\Rs$-algebras $\As\xrightarrow{\phi}\Bs$, there exists a unique covariant $\Rs$-linear homomorphism $\As^\gamma\xrightarrow{\hat{\phi}}\Bs$ such that $\phi=\hat{\phi}\circ\eta_\As$. 

\medskip 

In the case of $\Zs$-central unital $\Rs$-bimodules ${}_\Rs\Ms_\Rs$ the definition of $\gamma$-conjugate $\Ms\xrightarrow{\eta_\Ms}\Ms^\gamma$ is given via the same universal factorization property diagram of $\Rs$-bimodules, ``forgetting'' the multiplication. 
\end{definition}

\begin{remark}\label{rem: cj-dual}
Unicity of $\gamma$-conjugates up to a unique isomorphism compatible with the universal property is standard, their existence can be provided as follows. 

\medskip 

Given a $\Zs$-central unital associative $\Rs$-algebra $\As_\Rs$ and a conjugation $\gamma$ in $\Rs$, take as a $\Zs$-central bimodule $\As^\gamma:=\As$ and define $\eta_\As:\As\to\As^\gamma$ as the identity map, here denoted as $\As\ni x\mapsto \hat{x}\in\As^\gamma$. If $\gamma$ is a contravariant conjugation, define $r_1 \ \hat{\cdot}\ \hat{x} \ \hat{\cdot}\  r_2:=\widehat{\gamma(r_2)\cdot x\cdot \gamma(r_1)}$ and $\hat{x}\ \hat{\circ}\ \hat{y}:=\widehat{y\circ x}$, for all $x,y\in\As$ and $r_1,r_2\in\Rs$. If $\gamma$ is a covariant conjugation, define $r_1 \ \hat{\cdot}\ \hat{x} \ \hat{\cdot}\  r_2:=\widehat{\gamma(r_1)\cdot x\cdot \gamma(r_2)}$ and $\hat{x}\ \hat{\circ}\ \hat{y}:=\widehat{x\circ y}$, for all $x,y\in\As$ and $r_1,r_2\in\Rs$.

Notice that in both cases $\As^\gamma$ becomes a $\Zs$-central $\Rs$-bimodule with the new actions $\hat{\cdot}$ and it becomes a $\Zs$-central $\Rs$-algebra with the new product $\hat{\circ}$; furthermore the map $\eta_\As:\As\to\As^\gamma$ turns out to be a $\Zs$-linear $\gamma$-conjugate-linear contravariant (respectively covariant) unital homomorphism.  

For any $\gamma$-conjugate-linear unital contravariant (respectively contravariant) homomorphism $\phi:\As\to\Bs$, we necessarily need to define $\hat{\phi}(\hat{x}):=\phi(x)$, and we verify that $\hat{\phi}:\As^\gamma\to\Bs$ is an $\Rs$-linear unital covariant homomorphism in both cases. 
\xqed{\lrcorner}
\end{remark}

\section{Multimodules Over Unital Associative $\Zs$-central $\Rs$-algebras} \label{sec: multimodules}

We introduce here multimodules \textit{over families of unital associative $\Zs$-central $\Rs$-algebras}.\footnote{This generalizes the special case of multimodules over unital associative $\KK$-algebras over the field $\KK$: in this case one can take $\Rs:=\KK$ and $\Zs$ a subfield of $\KK$ consisting of fixed points for all the relevant conjugations of $\KK$ (in practice it is always possible, in each characteristic $p$, to take $\Zs$ as the initial field of that characteristic: $\QQ$ in characteristic $0$ and $\FF_p$ for any $p$ prime).}

\medskip 

In the following, we adapt the general definition of multimodule from~\cite[section~II.1.14]{Bou89}: 
\begin{definition}\label{def: morphism}
Let $\Zs$ be a commutative unital associative ring and $\Rs$ be a unital associative $\Zs$-central ring.
Given two families of unital associative $\Zs$-central $\Rs$-algebras $(\As_\alpha)_{\alpha\in A}$ and $(\Bs_\beta)_{\beta\in B}$, an \emph{$(\As_\alpha)$-$(\Bs_\beta)$ multimodule} ${}_{(\As_\alpha)}\Ms_{(\Bs_\beta)}$ is a $\Zs$-central bimodule that is a $\Zs$-central unital $\As_\alpha$-$\Bs_\beta$ bimodule for every $(\alpha,\beta)\in A\times B$ such that every pair of left actions and every pair of right actions commute.\footnote{
We assume the existence of a \textit{common} $\Zs$-central bimodule structure on $\Ms$ compatible with all the $\Zs$-bilinear right/left actions. 
} 

\medskip 

A \emph{morphism of multimodules} ${}_{(\As_\alpha)_A}\Ms_{(\Bs_\beta)_B}\xrightarrow{(\phi,\eta,\Phi,\zeta,\psi)_f}{}_{(\Cs_\gamma)_C}\Ns_{(\Ds_\delta)_D}$, $(A_+,B_+)$-covariant in the sub-families of indexes $A_+\subset A$,  $B_+\subset B$ and $(A_-,B_-)$-contravariant in the sub-families of indexes $A_-:=A-A_+$, $B_-:=B-B_+$, consists of: 
\begin{itemize}
\item
an injective function $f:A\uplus B\to C\uplus D$, with $A_+=A\cap f^{-1}(C)$, $B_+=B\cap f^{-1}(D)$; 
\item 
two maps $A\xrightarrow{\phi}\End_\Zs(\Rs)\xleftarrow{\psi}B$
associating to every pair of indexes $\alpha\in A$ and $\beta\in B$ two $\Zs$-linear unital endomorphisms $\phi_\alpha,\psi_\beta$ of $\Rs_\Zs$, covariant for $(\alpha,\beta)\in A_+\times B_+$ and contravariant for $(\alpha,\beta)\in A_-\times B_-$, 
\item 
for $(\alpha,\beta)\in A_+\times B_+$, 
$\Zs$-linear covariant unital homomorphisms $\As_\alpha\xrightarrow{(\phi_\alpha,\eta_\alpha)}\Cs_{f(\alpha)}$, $\Bs_\beta\xrightarrow{(\psi_\beta,\zeta_\beta)}\Ds_{f(\beta)}$; 
\item
for $(\alpha,\beta)\in A_-\times B_-$, 
$\Zs$-linear contravariant unital homomorphisms $\As_\alpha\xrightarrow{(\phi_\alpha,\eta_\alpha)}\Ds_{f(\alpha)}$, $\Bs_\beta\xrightarrow{(\psi_\beta,\zeta_\beta)}\Cs_{f(\beta)}$;
\item  
a $\Zs$-linear map $\Ms\xrightarrow{\Phi}\Ns$ such that $\Phi(a\cdot x\cdot b)=\eta_\alpha(a)\cdot \Phi(x)\cdot \zeta_\beta(b)$, for all $(\alpha,\beta)\in A_+\times B_+$ and $\Phi(a\cdot x\cdot b)=\zeta_\beta(b)\cdot \Phi(x)\cdot \eta_\alpha(a)$, for all $(\alpha,\beta)\in A_-\times B_-$, $(a,b)\in\As\times\Bs$ and $x\in\Ms$. 
\end{itemize}

\medskip 
The \emph{signature} of the morphism is $(\phi,\eta,\zeta,\psi)_f$. The function $f$ is the \emph{covariance of signature} of the morphism and \emph{covariant morphisms} are those for which $f(A)\subset C$ and $f(B)\subset D$. The pair $(\phi,\psi)$ is the \emph{$\Rs$-linearity of the signature} of the morphism and \emph{$\Rs$-linear morphisms} are those for which both $\phi$ and $\psi$ are constant equal to $\id_\Rs$. In some cases we will denote by $\Phi^\sigma$ a morphism $(\phi,\eta,\Phi,\zeta,\psi)_f$with signature $\sigma=(\phi,\eta,\zeta,\psi)_f$. 

\medskip 

The \emph{composition of morphisms} 
${}_{(\As_\alpha)_A}\Ms_{(\Bs_\beta)_B}\xrightarrow{(\phi_2,\eta_2,\Phi_2,\zeta_2,\psi_2)_{f_2}}{}_{(\As'_{\alpha'})_{A'}}\Ns_{(\Bs'_{\beta'})_{B'}}
\xrightarrow{(\phi_1,\eta_1,\Phi_1,\zeta_1,\psi_1)_{f_1}} {}_{(\As''_{\alpha''})_{A''}}\Ps_{(\Bs''_{\beta''})_{B''}}$
of multimodules is given componentwise: 
\begin{equation*}
(\phi_1,\eta_1,\Phi_1,\zeta_1,\psi_1)_{f_1}\circ (\phi_2,\eta_2,\Phi_2,\zeta_2,\psi_2)_{f_2}:=
(\phi_1\circ\phi_2,\eta_1\circ\eta_2,\Phi_1\circ\Phi_2,\zeta_1\circ\zeta_2,\psi_1\circ\psi_2)_{f_1\circ f_2}.
\end{equation*} 
The \emph{identity} of a multimodule ${}_{(\As_\alpha)_A}\Ms_{(\Bs_\beta)_B}$ is the morphism 
$(\id_\Rs, (\id_{\As_\alpha})_A, \id_\Ms,(\id_{\Bs_\beta})_B,\id_\Rs)_{\id_{A\uplus B}}$.
\end{definition}

\begin{remark}\label{rem: z2}
The map $\Phi:\Ms\to\Ns$ between multimodules does not have an intrinsic covariance: for every left index $\alpha\in A$ and for every right index $\beta\in B$ the morphisms $(\eta_\alpha,\Phi)$ and $(\zeta_\beta,\Phi)$ are covariant or contravariant depending on the sign $\pm$ indicated in the subsets $A_\pm$ and $B_\pm$. 

\medskip 

Similarly, the map $\Phi:\Ms\to\Ns$ between multimodules is always $\Zs$-linear, but it does not have an intrinsic $\phi$-linearity with respect to $\Rs$ for a fixed $\Zs$-linear morphism $\phi$: for every left index $\alpha\in A$ and right index $\beta\in B$, the morphism $(\eta_\alpha,\Phi)$ is $\phi_\alpha$-linear and the morphism $(\zeta_\beta,\Phi)$ is $\psi_\beta$-linear. 

\medskip 

We have a category $\Mf_{[\Rs_\Zs]}$ of morphisms of $\Zs$-central multimodules over $\Rs_\Zs$-algebras with composition of morphisms defined componentwise. 
The subcategories of $\Mf_{[\Rs_\Zs]}$ consisting of $(\As_\alpha)_{\alpha\in A}$-$(\Bs_\beta)_{\beta\in B}$ multimodules, over the same two families of unital associative $\Rs_\Zs$-algebras, and morphisms given by $(\phi,\eta,\Phi,\zeta,\psi)_f$, with $f:=\id_{A\uplus B}$, 
$\phi_\alpha:=\id_\Rs=:\psi_\beta$ and $\eta_\alpha=\id_{\As_\alpha}$, $\zeta_\beta=\id_{\Bs_\beta}$ for all $(\alpha,\beta)\in A\times B$, are denoted by ${}_{(\As_\alpha)_A}\Mf_{(\Bs_\beta)_B}$. In case of $\As$-bimodules, we use the notation ${}_\As\Mf_\As$.
\xqed{\lrcorner}
\end{remark}

This essential remark explains why the study of multimodules cannot be ``reduced'' to the theory of bimodules. 
\begin{remark}\label{rem: hyper}
If $\Rs=\Zs$, it is common to dismiss the usage of multimodules ${}_{(\As_\alpha)}\Ms_{(\Bs_\beta)}$ in favor of their ``equivalent'' description as bimodules ${}_{\bigotimes^\Rs_\alpha\As_\alpha}\Ms_{\bigotimes^\Rs_\beta\Bs_\beta}$ over tensor product $\Rs$-algebras $\bigotimes^\Rs_{\alpha\in A}\As_\alpha$ and ${\bigotimes^\Rs_{\beta\in B}\Bs_\beta}$ since: 
\begin{quote}
if $\Rs=\Zs$, there is a categorical isomorphism between the sub-category $_{(\As_\alpha)_A}\Mf_{(\Bs_\beta)_B}$ of \textit{covariant $\Rs$-linear} morphisms of $(\As_\alpha)_A$-$(\Bs_\beta)_B$-multimodules and the category ${}_{\bigotimes^\Rs_\beta\As_\alpha}\Mf_{\bigotimes^\Rs_\beta\Bs_\beta}$ of \textit{covariant $\Rs$-linear} morphisms of bimodules over the $\Rs$-balanced tensor product of the $\Rs_\Zs$-algebras. 
\end{quote}
As soon as one considers morphisms of multimodules with arbitrary covariance $f$, it is actually impossible to impose a unique unital associative product on the $\Rs$-tensor product algebras in order to obtain a similar equivalent treatment via categories of bimodules. 

\medskip 

A perfectly possile alternative (that we do not pursue here) would be to work with the category of ``bimodules'' over \textit{hyper-$\Zs$-central $\Rs$-algebras}: $\Zs$-central bimodules $\bigotimes^\Zs_{\alpha\in A}\As_\alpha$ equipped with many different $\Rs$-actions (on each of the tensor-factors) and different $\Zs$-bilinear associative unital binary product operations suitably compatible with the $\Rs$-actions (see for example~\cite[section~5.3]{BCLS20}); but in this case multimodules need anyway to be used in order to define hyper-algebras. 
\xqed{\lrcorner}
\end{remark}

To a certain extent, the usage of general morphisms of multimodules (with arbitrary conjugation and convariance signatures as in definition~\ref{def: morphism}) can be avoided, replacing the target multimodule with a suitable ``twisted version'' (depending on the signatures of the original morphism) and obtaining as a result an $\Rs$-linear covariant morphism into such ``twisted multimodule''. The construction follows similar steps as in definition~\ref{def: cj-dual} and remark~\ref{rem: cj-dual} and it simultaneously extends to multimodules the notions of \textit{conjugate-dual}, \textit{opposite}, \textit{restriction of rings}, \textit{pull-back}.

\begin{definition}\label{def: twisted}
Let ${}_{(\As_\alpha)_A}\Ms_{(\Bs_\beta)_B}$ and ${}_{(\Cs_\gamma)_C}\Ns_{(\Ds_\delta)_D}$ be two $\Zs$-central multimodules over $\Rs_\Zs$-algebras, and let \hbox{$\sigma:=(\phi,\eta,\zeta,\phi)_f$} be a given signature for multimodule morphisms between $\Ms$ and $\Ns$. 

\medskip 

A \emph{$\sigma$-twisted multimodule}\footnote{
We might also write \textit{$\Phi$-twisted} of $\Ns$, for a morphism $\Ms\xrightarrow{\Phi}\Ns$, instead of $\sigma(\Phi)$-twisted, where $\sigma(\Phi)$ denotes the signature of $\Phi$. 
} of $\Ns$ consists of a morphism of multimodules ${}_{(\As_\alpha)_A}{\Ns^\sigma}_{(\Bs_\beta)_B}\xrightarrow{\Theta^\sigma_\Ns}{}_{(\Cs_\gamma)_C}\Ns_{(\Ds_\delta)_D}$, with signature $\sigma$, such that the following universal factorization property is satisfied: 
for any other morphism of multimodules  ${}_{(\As_\alpha)_A}{\Ms}_{(\Bs_\beta)_B}\xrightarrow{\Phi}{}_{(\Cs_\gamma)_C}\Ns_{(\Ds_\delta)_D}$, with signature $\sigma$, there exists a unique covariant $\Rs$-linear morphism of multimodules 
${}_{(\As_\alpha)_A}{\Ms}_{(\Bs_\beta)_B}\xrightarrow{\Phi^\sigma}{\Phi}_{(\As_\alpha)_A}{\Ns^\sigma}_{(\Bs_\beta)_B}$
in the category ${}_{(\As_\alpha)_A}{\Mf}_{(\Bs_\beta)_B}$ such that $\Phi=\Theta^\sigma_\Ns\circ \Phi^\sigma$.
\end{definition}

\begin{remark}\label{prop: cj-duals}
As any definition via universal factorizations, $\sigma$-twisted of a given multimodule are unique, up to a unique isomorphism compatible with the factorization property. A construction can be achieved as follows.
Consider $\Ns^\sigma:=\Ns$ as a $\Zs$-central bimodule and $\Theta_\Ns^\sigma:\Ns^\sigma\to\Ns$ as the identity map. For all $x\in \Ns$ we will denote by $x^\sigma\in\Ns^\sigma$ its corresponding element, hence $\Theta^\sigma_\Ns(x^\sigma)=x$, for all $x\in\Ns$. For all $\sigma$-covariant indexes $(\alpha_+,\beta_+)\in A_+\times B_+$, and $\sigma$-contravariant indexes $(\alpha_-,\beta_-)\in A_-\times B_-$, we define new actions on $\Ns^\sigma$: 
\begin{gather*}
a\cdot_{\alpha_+} x^\sigma\cdot_{\beta_+} b:= \left(\eta_{\alpha_+}(a)\cdot_{f(a_+)} x\cdot_{f(b_+)} \zeta_{b_+}(b)\right)^\sigma, 
\quad \forall (a,b)\in\As_{\alpha_+}\times\Bs_{\beta_+}, \ \forall x^\sigma\in\Ns^\sigma, 
\\
a\cdot_{\alpha_-} x^\sigma\cdot_{\beta_-} b:=
\left( \zeta_{b_-}(b)\cdot_{f(b_-)} x \cdot_{f(a_-)} \eta_{\alpha_-}(a)\right)^\sigma, 
\quad \forall (a,b)\in\As_{\alpha_-}\times\Bs_{\beta_-}, \ \forall x^\sigma\in\Ns^\sigma, 
\end{gather*}
obtaining a multimodule ${}_{(\As_\alpha)_A}{\Ns^\sigma}_{(\Bs_\beta)_B}$ such that the map $\Theta_\Ns^\sigma:x^\sigma\mapsto x$ is a morphism of multimodules with signature $\sigma$. 
Finally, given any other morphism ${}_{(\As_\alpha)_A}{\Ms}_{(\Bs_\beta)_B}\xrightarrow{\Phi}{}_{(\Cs_\gamma)_C}\Ns_{(\Ds_\delta)_D}$ of multimodules with signature $\sigma$, the function $\Phi^\sigma: m\mapsto (\Phi(m))^\sigma\in\Ns^\sigma$, (due to the bijectivity of $\Theta_\Ns^\sigma$) is the unique map that satisfies
$\Theta_\Ns^\sigma(\Phi^\sigma(m))=\Theta_\Ns^\sigma((\Phi(m))^\sigma)=\Phi(m)$, for all $m\in\Ms$, and by direct calculation,
we see that it is also a morphism of multimodules with identity signature. 
\xqed{\lrcorner}
\end{remark}

As typical of any category of homomorphisms of algebraic structures, sub-structures can be defined via algebraically closed subsets and quotient-structures via congruences. 
\begin{definition}
Given a multimodule ${}_{(\As_\alpha)_A}\Ms_{(\Bs_\beta)_B}$ over $\Zs$-central $\Rs$-algebras, 
\begin{itemize}
\item 
a \emph{sub-multimodule} of $\Ms$ is a subset $\Ns\subset\Ms$ that is \emph{algebraically closed} under all the operations: 
\begin{equation*}
0_\Ms\in\Ns, \quad 
x,y\in\Ns \imp x+y\in\Ns, 
\quad 
x\in\Ns \imp a\cdot_\alpha x\cdot_\beta b\in\Ns, 
\quad 
\forall (\alpha,\beta)\in A\times B, \ (a,b)\in\As_\alpha\times\Bs_\beta, \ x,y\in \Ns;
\end{equation*}
\item
a \emph{multimodule congruence} on ${}_{(\As_\alpha)_A}\Ms_{(\Bs_\beta)_B}$ is an equivalence relation $\E\subset\Ms\times\Ms$ such that: 
\begin{equation*}
x\sim_\E y \imp (x+z)\sim_\E (y+z), 
\quad 
x\sim_\E y \imp (a\cdot_\alpha x\cdot_\beta b)\sim_\E (a\cdot_\alpha y\cdot_\beta b), 
\quad  \forall (\alpha,\beta)\in A\times B, \ (a,b)\in\As_\alpha\times\Bs_\beta, \ x,y,z\in \Ms. 
\end{equation*}
\end{itemize}
A \emph{quotient multimodule} of $\Ms$ by the congruence $\E$ is the multimodule ${}_{(\As_\alpha)_A}(\frac{\Ms}{\E})_{(\Bs_\beta)_B}$ consisting of the quotient set $\Ms/\E$ equipped with the well-defined addition $[x]_\E+[y]_\E:=[x+y]_\E$, for all $x,y\in\Ms$, and the well-defined actions: 
$a\cdot_\alpha [x]_\E\cdot_\beta b:=[a\cdot_\alpha x\cdot_\beta b]_\E$, for all $x\in\Ms$, $(\alpha,\beta)\in A\times B$ and $(a,b)\in\As_\alpha\times\Bs_\beta$.  
\end{definition}
\begin{remark}
As usual, any multimodule congruence $\E$ uniquely determines the $\Ms$-sub-multimodule $[0_\Ms]_\E$; reciprocally any $\Ms$-sub-multimodule $\Ns$ uniquely determines a multimodule congruence $x\sim y :\iff x-y\in\Ns$ whose equivalence classes, for $x\in\Ms$, are the affine spaces $[x]_\sim=x+\Ns:=\{x+y \ | \ y\in\Ns \}$. The notation $\frac{\Ms}{\Ns}$ is used to identify the quotient of $\Ms$ by the congruence uniquely determined by the sub-multimodule $\Ns$. 

\medskip 

Inclusions of sub-multimodules $\Ns\xrightarrow{\iota}\Ms$ and quotients $\Ms\xrightarrow{\pi}\frac{\Ms}{\Ns}$ are morphisms in the category ${}_{(\As_\alpha)_A}\Mf_{(\Bs_\beta)_B}$. 
\xqed{\lrcorner}
\end{remark}

\medskip 

Despite being rarely mentioned, multimodules naturally appear whenever bimodules are around: 

\begin{proposition}\label{prop: homK}
Let ${}_\As\Ms_\Bs$ and ${}_{\As'}\Ns_{\Bs'}$ be $\Zs$-central bimodules over $\Zs$-central unital associative $\Rs$-algebras $\As, \As',\Bs,\Bs'$. The set $\Hom_\Zs(\Ms;\Ns)$ of \hbox{$\Zs$-linear} maps $\phi:\Ms\to\Ns$ is a left-$(\As',\Bs)$ right-$(\As,\Bs')$ multimodule with the following actions, for all $a\in\As$, $a'\in\As'$, $x\in\Ms$, $b\in\Bs$, $b'\in\Bs'$: 
\begin{align*}
&\textit{left external action:} &&(a'\cdot \phi)(x):=a'\phi(x), 
\\
&\textit{right external action:} &&(\phi\cdot b')(x):=\phi(x)b', 
\\
&\textit{left internal action:} &&(b\odot \phi)(x):=\phi(xb), 
\\
&\textit{right internal action:} &&(\phi\odot a)(x):=\phi(ax). 
\end{align*}
\end{proposition}
\begin{proof}
By direct calculation, $x\mapsto (a'\cdot\phi)(x)$, $x\mapsto (\phi\cdot b')(x)$, $x\mapsto (b\odot \phi)(x)$, $x\mapsto(\phi\odot a)(x)$ are all $\Zs$-linear and the above defined maps are all $\Zs$-bilinear actions. 
To prove the multimodule structure on $\Hom_\Zs(\Ms;\Ns)$ we check that the actions pairwise commute, for all $a\in\As$, $a'\in\As'$ and all $b\in\Bs$, $b'\in\Bs'$:
\begin{align*}
&(a'\cdot\phi)\cdot b'=a'\cdot(\phi\cdot b'), 
&&
(a'\cdot\phi)\odot a=a'\cdot(\phi\odot a), 
\\
&(b\odot\phi)\odot a=b\odot(\phi\odot a), 
&& 
(b\odot\phi)\cdot b'=b\odot(\phi\cdot b'). 
\qedhere 
\end{align*}
\end{proof}

\begin{remark}\label{rem: Hom-multi}
More generally, if ${}_{(\As_\alpha)}\Ms_{(\Bs_\beta)}$ and ${}_{(\Cs_\gamma)}\Ns_{(\Ds_\delta)}$ are $\Zs$-central multimodules over $\Zs$-central $\Rs$-algebras, the $\Zs$-central bimodule $\Hom_\Zs(\Ms;\Ns)$ becomes a left-$(\Bs_\beta,\Cs_\gamma)_{\beta\in B,\gamma\in C}$ and a right-$(\As_\alpha,\Ds_\delta)_{\alpha\in A,\delta\in D}$ $\Zs$-central multimodule with internal/external actions given by: 
\begin{gather*}
(c_\gamma\cdot \phi\cdot d_\delta)(x):=c_\gamma\cdot \phi(x)\cdot d_\delta, 
\quad 
(b_\beta\odot \phi\odot a_\alpha)(x):=\phi(a_\alpha\cdot x\cdot b_\beta), 
\end{gather*}
for all $(\alpha,\beta,\gamma,\delta)\in A\times B\times C\times D$,  $(a_\alpha,b_\beta,c_\gamma,d_\delta)\in\As_\alpha\times\Bs_\beta\times\Cs_\gamma\times\Ds_\delta$, $x\in\Ms$ and $\phi\in\Hom_\Zs(\Ms;\Ns)$. 
\xqed{\lrcorner}
\end{remark}

Tensor products provide other examples of multimodules~\cite[section~II.3.4]{Bou89}: 
\begin{proposition}
Let ${}_{(\As_\alpha)}\Ms_{(\Bs_\beta)}$ and ${}_{(\Cs_\gamma)}\Ns_{(\Ds_\delta)}$ be $\Zs$-central multimodules over $\Zs$-central $\Rs$-algebras. Their tensor product $\Ms\otimes_\Zs\Ns$ over $\Zs$ is a left-$(\As_\alpha,\Cs_\gamma)_{\alpha\in A,\gamma\in C}$ right-$(\Bs_\beta,\Ds_\delta)_{\beta\in B,\delta\in D}$ $\Zs$-central multimodule.
\end{proposition}
\begin{proof}
The definition of tensor product (via universal factorization property for $\Zs$-balanced bi-homomorphism) and its construction are well-known: see for example~\cite[section~II.3, proposition~3]{Bou89}; we only recall here the relevant actions on simple tensors: 
\begin{equation*}
a\cdot (x\otimes_\Zs y)=(a\cdot x)\otimes_\Zs y, 
\quad 
c\cdot (x\otimes_\Zs y)=x\otimes_\Zs (c\cdot y), 
\quad 
(x\otimes_\Zs y)\cdot b=(x\cdot b)\otimes_\Zs y, 
\quad 
(x\otimes_\Zs y)\cdot d=x\otimes_\Zs (y\cdot d),  
\end{equation*}
for all $(x,y)\in\Ms\times\Ns$, $(a,b,c,d)\in\As_\alpha\times\Bs_\beta\times\Cs_\gamma\times\Ds_\delta$, $(\alpha,\beta,\gamma,\delta)\in A\times B\times C\times D$. 
\end{proof}

\begin{remark}\label{rem: otimes-tr}
One can actually define tensor products of multimodules in much greater generality. 

\medskip 

Instead of taking only the tensor product over the algebra $\Zs$ of ``scalars'' and use $\Zs$-bilinear maps, we can ``contract'' over arbitrary families of shared $\Zs$-central $\Rs$-algebras acting on the two multimodules and utilize suitable maps that are ``balanced'' over the ``contracted actions'', obtaining multimodules over the remaining ``un-contracted'' actions, as detailed in the following exposition. 

\medskip 

Let ${}_{(\As_\alpha)_A}\Ms_{(\As_\beta)_B}$ and ${}_{(\As_\gamma)_C}\Ns_{(\As_\delta)_D}$ be a pair of multimodules; consider the relation $A\uplus B\xrightarrow{\Sigma}C\uplus D$, defined by $(\xi,\zeta)\in\Sigma \iff \As_\xi=\As_\zeta$ (where $\xi\in A\uplus B$ and $\zeta\in C\uplus D$), and let $A\uplus B \supset A'\uplus B' \xrightarrow{\Gamma} C'\uplus D'\subset C\uplus D$ be a bijective function between subsets of indexes, such that $\Gamma\subset \Sigma$ (in practice $\As_\xi=\As_{\Gamma(\xi)}$, for all $\xi\in A'\uplus B'$). 

\medskip 

A \emph{tensor product of multimodules} $\Ms$ and $\Ns$ over $\Gamma$ consists of: 
\begin{itemize}
\item 
a left-$(\As_\xi)_{\xi\in (A-A')\uplus(C-C')}$ right-$(\As_\zeta)_{\zeta\in (B-B')\uplus(D-D')}$ multimodule 
${}_{(\As_\xi)_{(A-A')\uplus(C-C')}}(\Ms\otimes_\Gamma\Ns)_{(\As_\zeta)_{(B-B')\uplus(D-D')}}$, 
\item 
a \emph{$\Gamma$-balanced bi-morphism} $\Ms\times\Ns\xrightarrow{\eta} \Ms\otimes_\Gamma\Ns$, that means a $\Zs$-bilinear map that satisfies: 
\begin{align} \notag
&
\eta(a \cdot_\xi x,y)=\eta(x,a\cdot_{\Gamma(\xi)} y), \ \forall a\in\As_\xi, \ (\xi,\Gamma(\xi))\in A'\times C', \ (x,y)\in\Ms\times\Ns, 
\\ \notag
&
\eta(x \cdot_\xi a,y)=\eta(x,y\cdot_{\Gamma(\xi)} a), \ \forall a\in\As_\xi, \ (\xi,\Gamma(\xi))\in B'\times D', \ (x,y)\in\Ms\times\Ns, 
\\ \notag
&\eta(x \cdot_\xi a,y)=\eta(x,a\cdot_{\Gamma(\xi)} y), \ \forall a\in\As_\xi, \ (\xi,\Gamma(\xi))\in B'\times C', \ (x,y)\in\Ms\times\Ns, 
\\ \label{eq: balanced-cong}
&\eta(a \cdot_\xi x,y)=\eta(x,y\cdot_{\Gamma(\xi)} a), \ \forall a\in\As_\xi, \ (\xi,\Gamma(\xi))\in A'\times D', \ (x,y)\in\Ms\times\Ns, 
\\ \notag
&\eta(a \cdot_\xi x,y)=a\cdot_\xi\eta(x,y), \ \forall a\in\As_\xi, \ \xi\in A-A', \ (x,y)\in\Ms\times\Ns,
\\ \notag
&\eta(x,c \cdot_\xi y)=c\cdot_\xi\eta(x,y), \ \forall c\in\As_\xi, \ \xi\in C-C', \ (x,y)\in\Ms\times\Ns,
\\ \notag
&\eta(x \cdot_\xi b,y)=\eta(x,y)\cdot_\xi b, \ \forall b\in\As_\xi, \ \xi\in B-B', \ (x,y)\in\Ms\times\Ns,
\\ \notag
&\eta(x,y \cdot_\xi d)=\eta(x,y)\cdot_\xi d, \ \forall d\in\As_\xi, \ \xi\in D-D', \ (x,y)\in\Ms\times\Ns,
\end{align} 
\end{itemize}
in such a way that the following universal factorization property holds: for any other $\Gamma$-balanced bi-morphism  
$\Ms\times\Ns\xrightarrow{\Phi}\Ps$ into a left-$(\As_\xi)_{(A-A')\uplus(C-C')}$ right-$(\As_\zeta)_{\zeta\in (B-B')\uplus(D-D')}$ multimodule $\Ps$, there exists a unique morphism of multimodules $\Ms\otimes_\Gamma\Ns\xrightarrow{\hat{\Phi}}\Ps$ (over the same indexed families of algebras) such that $\Phi=\hat{\Phi}\circ \eta$. 

\medskip 

Its construction is standard and consists of the quotient of a free multimodule over $\Ms\times\Ns$ by the congruence generated by the required axioms of $\Gamma$-balanced bi-morphism. More specifically we recall that: 
\begin{itemize} 
\item 
a \emph{free $(\As_\alpha)_{\alpha\in A}$-$(\Bs_\beta)_{\beta\in B}$ multimodule}, over a set $X$, is function  $X\xrightarrow{\eta^X}{}_{(\As_\alpha)_A}\Fg(X)_{(\Bs_\beta)_B}$, with values into a  
$(\As_\alpha)_{\alpha\in A}$-$(\Bs_\beta)_{\beta\in B}$ multimodule $\Fg(X)$, such that the following universal factorization property is satisfied: 
for any other map $X\xrightarrow{\Phi}{}_{(\As_\alpha)_A}\Ms_{(\Bs_\beta)_B}$ into an $(\As_\alpha)_{\alpha\in A}$-$(\Bs_\beta)_{\beta\in B}$ multimodule $\Ms$, there exists a unique morphism of multimodules $\Fg(X)\xrightarrow{\hat{\Phi}}\Ms$ in the category ${}_{(\As_\alpha)_A}\Mf_{(\Bs_\beta)_B}$ such that $\phi=\hat{\Phi}\circ\eta^X$; 
\item 
a construction of free multimodule over $X$ can be achieved taking $\bigoplus_{x\in X}[(\otimes^\Zs_{\alpha\in A}\As_\alpha)\otimes_\Zs(\otimes^\Zs_{\beta\in B}\Bs)]$, the set of finitely supported functions from $X$ into 
$(\otimes^\Zs_{\alpha\in A}\As_\alpha)\otimes_\Zs(\otimes^\Zs_{\beta\in B}\Bs)$ with pointwise addition and pointwise outer target actions as specified in footnote~\ref{foo: free-m}, defining $\eta^X_x(y):=\begin{cases}
(\otimes^\Zs_{\alpha\in A}1_{\As_\alpha})\otimes_\Zs (\otimes^\Zs_{\beta\in B}1_{\Bs_\beta}), \ y=x, 
\\ 
(\otimes^\Zs_{\alpha\in A}0_{\As_\alpha})\otimes_\Zs (\otimes^\Zs_{\beta\in B}0_{\Bs_\beta}), \ y\neq x,
\end{cases}$ for all $x,y\in X$, and checking the universal factorization property; 
\item 
the congruence $\E_\Gamma$ generated by the relations in equations~\eqref{eq: balanced-cong} is just the intersection of the set of congruences of $(\As_\alpha)_A$-$(\Bs_\beta)_B$ multimodule in $\Fg(\Ms\times \Ns)$, that contain all of the differences between left and right terms in each of the equations~\ref{eq: balanced-cong}; 
\item 
the tensor product consists of the quotient multimodule $\Ms\otimes_\Gamma\Ns:=\frac{\Fg(\Ms\times\Ns)}{\E_\Gamma}$, with the $\Gamma$-balanced bi-morphism $\eta:=\pi\circ\eta^{\Ms\times\Ns}$, where $\Fg(\Ms\times\Ns)\xrightarrow{\pi}\Ms\otimes_\Gamma\Ns$ is the quotient morphism. 
\xqedhere{3cm}{\lrcorner}
\end{itemize}
\end{remark}

\section{Involutions in Multimodules} \label{sec: involutions-m}

In parallel with the case of morphisms, also the nature of involutions in multimodules is more delicate and 
an involution is an involutive endomorphism inducing involutions on the algebras and conjugations on $\Rs$. 
\begin{definition}
Let ${}_{(\As_\alpha)_A}\Ms_{(\As_\beta)_B}$ be a $\Zs$-central multimodule over $\Zs$-central unital associative $\Rs$-algebras. 

\medskip

A \emph{multimodule involution} on $\Ms$ is a morphism  ${}_{(\As_\alpha)_A}\Ms_{(\As_\beta)_B}\xrightarrow{(\phi,\eta,\star,\zeta,\psi)_f}{}_{(\As_\alpha)_A}\Ms_{(\As_\beta)_B}$  that is involutive: 
\begin{itemize}
\item 
$A\uplus B\xrightarrow{f}A\uplus B$ is an involutive function $f\circ f=\id_{A\uplus B}$;\footnote{
From the involutivity of $f$, we have $f(A_+)=A_+, \ f(B_+)=B_+, \ f(A_-)=B_-$ and $f(B_-)=A_-$. 
} 
\item 
for all $(\alpha_1,\alpha_2)\in f\cap (A\times A)$, $\As_{\alpha_1}=\As_{\alpha_2}$, 
$\dagger_{\alpha_1}:=\phi_{\alpha_1}=\phi_{\alpha_2}$ is a covariant $\Rs_\Zs$-conjugation, 
$\ddagger_{\alpha_1}:=\eta_{\alpha_1}=\eta_{\alpha_2}$ is a covariant 
$\dagger_{\alpha_1}$-linear involution; 
\item 
for all $(\beta_1,\beta_2)\in f\cap (B\times B)$, $\As_{\beta_1}=\As_{\beta_2}$, 
$\dagger_{\beta_1}:=\psi_{\beta_1}=\psi_{\beta_2}$ is a covariant $\Rs_\Zs$-conjugation, 
$\ddagger_{\beta_1}:=\zeta_{\beta_1}=\zeta_{\beta_2}$ is a covariant 
$\dagger_{\beta_1}$-linear involution; 
\item
for all $(\alpha,\beta)\in f\cap (A\times B)$, $\As_{\alpha}=\As_{\beta}$, 
$\dagger_\alpha:=\phi_{\alpha}=\psi_{\beta}=:\dagger_\beta$ is a contravariant $\Rs_\Zs$-conjugation, 
$\ddagger_\alpha:=\eta_{\alpha}=\zeta_{\beta}=:\ddagger_\beta$ is a contravariant 
$\dagger_\alpha$-linear involution; 
\item 
${}_{\As_{\alpha}}\Ms_{\As_{\beta}}\xrightarrow{\star}{}_{\As_{\alpha}}\Ms_{\As_{\beta}}$ is an involution such that:
\begin{gather*}
\forall (\alpha,\beta)\in A_+\times B_+ \st (a\cdot_\alpha x\cdot_\beta b)^\star= a^{\ddagger_\alpha}\cdot_{f(\alpha)} x^\star\cdot_{f(\beta)} b^{\ddagger_\beta}, \quad \forall (a,b)\in \As_\alpha\times\As_\beta, \ x\in\Ms; 
\\
\forall (\alpha,\beta)\in A_-\times B_- 
\st (a\cdot_\alpha x\cdot_\beta b)^\star=b^{\ddagger_\beta}\cdot_{f(\beta)} x^\star\cdot_{f(\alpha)} a^{\ddagger_\alpha}, \quad \forall (a,b)\in\As_\alpha\times\As_\beta, \ x\in\Ms. 
\end{gather*} 
\end{itemize}
If necessary, we will denote an involution of ${}_{(\As_\alpha)_{A_\Ms}}\Ms_{(\As_\beta)_{B_\Ms}}$ by 
$(\dagger^\sigma_\Ms,\ddagger^\sigma_\Ms,\star_\Ms)_{\sigma\in f_\Ms}$, where: 
\begin{equation*}
\vcenter{\xymatrix{
\Zs \ar@{=}[d] \ar[r]^{\iota_\Rs} & \Rs \ar[d]^{\dagger_\Ms^\sigma} \ar[r]^{\iota_\As} & \As_{\sigma_1} \ar[d]^{\ddagger_\Ms^\sigma} & \Ms \ar[d]^{\star_\Ms}
\\
\Zs \ar[r]^{\iota_\Rs} & \Rs \ar[r]^{\iota_\Bs} & \As_{\sigma_2}  & \Ms 
}}\quad \sigma:=(\sigma_1,\sigma_2)\in f_\Ms\subset (A_\Ms\uplus B_\Ms)\times(A_\Ms\uplus B_\Ms). 
\end{equation*}  
\end{definition} 

Here we examine involutions for multimodules of morphisms between involutive bimodules. 
\begin{proposition}\label{prop: dagger}
Suppose that $\As$ and $\Bs$ are both $\Zs$-central $\Rs_\Zs$-algebras with involutions $\ddagger_\As$ and $\ddagger_\Bs$ over the respective $\Rs$-conjugations $\dagger_\As$ and $\dagger_\Bs$. 
If the $\Zs$-central bimodules ${}_\As\Ms_\As$ and ${}_{\Bs}\Ns_{\Bs}$ are both involutive with involutions $(\dagger_\As,\ddagger_\As,\star_\Ms)_{f_\Ms}$ and $(\dagger_\Bs,\ddagger_\Bs,\star_\Ns)_{f_\Ns}$,  
also the $\Zs$-central multimodule ${}_{\Bs,\As}\Hom_\Zs(\Ms;\Ns)_{\As,\Bs}$, considered in proposition~\ref{prop: homK}, is equipped with an involutive map $\star:T\mapsto T^\star:=\star_\Ns\circ T\circ \star_\Ms$ and becomes an involutive multimodule with multimodule involution $(\dagger,\ddagger,\star)_f$, defined as follows: 
\begin{equation*}
\vcenter{\xymatrix{
\Rs \ar[d]_{\dagger^\rho:=\dagger_\As} \ar[r]^{\iota_\As} & \As_{} \ar[d]^{\ddagger^\rho:=\ddagger_\As} & 
\Hom_\Zs(\Ms;\Ns)  \ar[d]^{\star}
& \Bs_{} \ar[d]_{\ddagger^\sigma:=\ddagger_\Bs} & \Rs \ar[d]^{\dagger^\sigma:=\dagger_\Bs} \ar[l]_{\iota_\Bs}
\\
\Rs \ar[r]^{\iota_\As} & \As_{}  & \Hom_\Zs(\Ms;\Ns) 
& \Bs_{} & \Rs  \ar[l]_{\iota_\Bs}
}}
\quad 
f:=f_\Ms^*\uplus f_\Ns, \quad \rho\in f_\Ms^*, \ \sigma\in f_\Ns. 
\end{equation*}

The involution $\star$ has covariance signature and $\Rs$-linearity signatures that, for inner actions, coincide with those of $\star_\Ms$; and for outer actions with those of $\star_\Ns$.

\medskip  

If $({}_{\Cs}\Ps_{\Cs},\star_\Ps)$ is an involutive $\Zs$-central bimodule over $\Rs_\Zs$-algebras, we have 
$(T\circ S)^\star=T^\star\circ S^\star$, for all $(T,S)\in \Hom_\KK(\Ns;\Ps)\times\Hom_\KK(\Ms;\Ns)$. In particular $(\Hom_\Zs(\Ms;\Ms),\circ,\star)$ is a unital associative $\Zs$-central algebra with a covariant involution.\footnote{
Notice that the involution $\star$ is multiplicative independently from the convariace/contravariace of the original involutions on $\Ms$.} 
\end{proposition}
\begin{proof}
If $T\in\Hom_\Zs(\Ms;\Ns)$ with $\Rs$-linearity signature $\phi_T$, the composition $T^\star:=\star_\Ns\circ T\circ \star_\Ms$ is $\Zs$-linear and with $\Rs$-linearity signature $\dagger_\Ns\circ\phi_T\circ\dagger_\Ms$ and hence 
$T\mapsto T^\star$ is well-defined as an endo-map of $\Hom_\Zs(\Ms;\Ns)$. 

\medskip 
 
The involutivity of $\star$ follows from: $(T^{\star})^{\star}=\star_\Ns\circ\star_\Ns\circ T\circ  \star_\Ms\circ \star_\Ms=T$, for all $T\in\Hom_\Zs(\Ms;\Ns)$.

\medskip 

For the actions, if $(\As,\ddagger_\As)$ and $(\Bs,\ddagger_\Bs)$ are contravariantly involutive, we necessarily have: 
\begin{align*}
(c\cdot b\odot T\odot a\cdot d)^{\star}(x)&=(c\cdot T(a\cdot x^{\star_\Ms}\cdot b)\cdot d)^{\star_\Ns}
=(d^{\ddagger_\Bs}\cdot T((b^{\ddagger_\As}\cdot x\cdot a^{\ddagger_\As})^{\star_\Ms})^{\star_\Ns}\cdot c^{\ddagger_\Bs})
\\
&=(d^{\ddagger_\Bs}\cdot a^{\ddagger_\As}\odot T^\star\odot b^{\ddagger_\As}\cdot c^{\ddagger_\Bs})(x), 
\quad \forall a,b\in\As, \ c,d\in \Bs, \ x\in\Ms. 
\end{align*}
Whenever $(\As,\ddagger_\As)$ and $(\Bs,\ddagger_\Bs)$ are covariantly involutive, we obtain: 
\begin{align*}
(c\cdot b\odot T\odot a\cdot d)^\star(x)&=(c\cdot T(a\cdot x^{\star_\Ms}\cdot b)\cdot d)^{\star_\Ns}
=(c^{\ddagger_\Bs}\cdot T((a^{\ddagger_\As}\cdot x\cdot b^{\ddagger_\As})^{\star_\Ms})^{\star_\Ns}\cdot d^{\ddagger_\Bs})
\\
&=(c^{\ddagger_\Bs}\cdot b^{\ddagger_\As}\odot T^\star\odot a^{\ddagger_\As}\cdot d^{\ddagger_\Bs})(x), 
\quad \forall a,b\in\As, \ c,d\in \Bs, \ x\in\Ms. 
\end{align*}
The remaining two cases with opposite contravariance between $(A,\ddagger_\As)$ and $(\Bs,\ddagger_\Bs)$ are treated similarly. 

\medskip 

Finally $(T\circ S)^\star=\star_\Ps\circ T\circ S\circ \star_\Ms =\star_\Ps\circ T\circ \star_\Ns\circ\star_\Ns \circ S\circ \star_\Ms=T^\star\circ S^\star, \ \forall (T,S)\in\Hom_\KK(\Ns;\Ps)\times\Hom_\KK(\Ms;\Ns)$.

\medskip 

Notice that the involution $\star:T\mapsto T^\star$ has $\Rs$-linearity signature $\dagger_\Ms$ for the inner actions and the $\Rs$-linearity signature of $\dagger_\Ns$ for the outer actions. 
\end{proof}

\begin{remark}
The previous proposition can easily be further generalized: whenever ${}_{(\As_\alpha)}\Ms_{(\As_\beta)}$ and ${}_{(\Bs_\gamma)}\Ns_{(\Bs_\delta)}$ are $\Zs$-central multimodules over $\Zs$-central $\Rs_\Zs$-algebras, any pair $(\dagger_\Ms,\ddagger_\Ms,\star_\Ms)_{f_\Ms}$ and $(\dagger_\Ns,\ddagger_\Ns,\star_\Ns)_{f_\Ns}$  of involutions, 
induces an involution $\star:T\mapsto T^\star:=\star_\Ns\circ T\circ\star_\Ms$ of $\Hom_\Zs(\Ms;\Ns)$, 
that is compatible with all the external and internal actions of the multimodule ${}_{(\Bs_\gamma,\As_\beta)}\Hom_\KK(\Ms;\Ns)_{(\As_\alpha,\Bs_\delta)}$ defined in remark~\ref{rem: Hom-multi} and hence, defining $f:=f_\Ms^*\uplus f_\Ns$, $\dagger:=\dagger_\Ms\uplus \dagger_\Ns$, $\ddagger:=\ddagger_\Ms\uplus\ddagger_\Ns$, we see that $(\dagger,\ddagger,\star)_f$ is an involution of the $\Zs$-central multimodule ${}_{(\Bs_\gamma,\As_\beta)}\Hom_\KK(\Ms;\Ns)_{(\As_\alpha,\Bs_\delta)}$ over $\Rs_\Zs$-algebras.\footnote{
Here, given two functions $F:A\to B$ and $G:C\to D$ with define $F\uplus G:A\uplus B\to C\uplus D$ the ``disjoint union'' of the two maps.
}  
\xqed{\lrcorner}
\end{remark}

The following proposition describes involutions in the case of tensor products of involutive multimodules. 
\begin{proposition}\label{prop: inv-ot}
Let $\left({}_{(\As_\alpha)_A}\Ms_{(\As_\beta)_B},((\star_\alpha)_A,\star_\Ms,(\star_\beta)_B)_f\right)$ and $\left({}_{(\Bs_\gamma)_C}\Ns_{(\Bs_\delta)_D},((\star_\gamma)_C,\star_\Ns,(\star_\delta)_D)_g\right)$ be involutive $\Zs$-central multimodules over $\Rs_\Zs$-algebras; 
the tensor product multimodule 
${}_{(\As_\alpha,\Bs_\gamma)_{A\uplus C}}(\Ms\otimes_\Zs\Ns)_{(\As_\beta,\Bs_\delta)_{B\uplus D}}$ has an involution $((\star_\alpha,\star_\beta)_{A\uplus C}, \star_\Ms\otimes_\Zs\star_\Ns,(\star_\gamma,\star_\delta)_{B\uplus D})_{(f,g)}$. 
\end{proposition}
\begin{proof} 
Define $A\uplus B\uplus C\uplus D \xrightarrow{(f,g)}A\uplus B\uplus C\uplus D$ as the ``disjoint union'' of the involutions $A\uplus B\xrightarrow{f} A\uplus B$ and $C\uplus D\xrightarrow{g} C\uplus D$. It follows that $(f,g)$ is an involution. 
Furthermore, for all $\Ts\in\{\As,\Bs\}$, for all $\tau\in\{\alpha,\beta,\gamma,\delta\}$ we have  
$(\Ts_{\tau},\star_{\tau})=(\Ts_{(f,g)(\tau)}\star_{(f,g)(\tau)})$. 

\medskip 

The $\Zs$-linear map $\Ms\otimes_\Zs\Ns\xrightarrow{\star_\Ms\otimes_\Zs\star_\Ns}\Ms\otimes_\Zs\Ns$, defined by universal factorization property from the $\Zs$-bilinear map 
$\Ms\times\Ns \ni (x,y)\mapsto x^{\star_\Ms}\otimes_\Zs y^{\star_\Ns}\in\Ms\otimes_\Zs\Ns$, for all $x\in\Ms$, $y\in\Ns$, is involutive. 

\medskip 

The covariance/contravariance behavior of the involution $\star_\Ms\otimes_\Zs\star_\Ns$ with respect to the several actions is described as follows, denoting by $\tau_\pm$, for $\tau\in\{\alpha,\beta,\gamma,\delta\}$, the indexes corresponding respectively to covariantly/contravariantly involutive algebras: 
\begin{align*}
\forall (\alpha_+,\beta_+,\gamma_+,\delta_+)\in A_+\times B_+\times C_+\times D_+, \ 
&\forall (a,b)\in\As_{\alpha_+}\times\As_{\beta_+},\ (c,d)\in\Bs_{\gamma_+}\times\Bs_{\delta_+},\ (x,y)\in\Ms\times\Ns:
\\
(a\cdot_{\alpha_+} c\cdot_{\gamma_+}(x\otimes_\Zs y) \cdot_{\beta_+} b \cdot_{\delta_+} d)^{\star_\Ms\otimes_\Zs\star_\Ns}
&=(a\cdot_{\alpha_+}x\cdot_{\beta_+} b)^{\star_\Ms}\otimes_\Zs (c\cdot_{\gamma_+} y \cdot_{\delta_+} d)^{\star_\Ns}
\\
&=(a^{\star_{\alpha_+}}\cdot_{f(\alpha_+)}x^{\star_\Ms}\cdot_{f(\beta_+)} b^{\star_{\beta_+}})\otimes_\Zs (c^{\star_{\gamma_+}}\cdot_{g(\gamma_+)} y^{\star_\Ns} \cdot_{g(\delta_+)} d^{\star_{\delta_+}})
\\
&=a^{\star_{\alpha_+}}\cdot_{f(\alpha_+)} c^{\star_{\gamma_+}}\cdot_{g(\gamma_+)}
(x\otimes_\Zs y)^{\star_\Ms\otimes_\Zs\star_\Ns} 
\cdot_{f(\beta_+)} b^{\star_{\beta_+}}\cdot_{g(\delta_+)} d^{\star_{\delta_+}}, 
\\
\forall (\alpha_-,\beta_-,\gamma_-,\delta_-)\in A_-\times B_-\times C_-\times D_-, \ 
&\forall (a,b)\in\As_{\alpha_-}\times\As_{\beta_-},\ (c,d)\in\Bs_{\gamma_-}\times\Bs_{\delta_-},\ (x,y)\in\Ms\times\Ns:
\\
(a\cdot_{\alpha_-} c\cdot_{\gamma_-}(x\otimes_\Zs y) \cdot_{\beta_-} b \cdot_{\delta_-} d)^{\star_\Ms\otimes_\Zs\star_\Ns}
&=(a\cdot_{\alpha_-}x\cdot_{\beta_-} b)^{\star_\Ms}\otimes_\Zs (c\cdot_{\gamma_-} y \cdot_{\delta_-} d)^{\star_\Ns}
\\
&=(b^{\star_{\beta_-}}\cdot_{f(\beta_-)}x^{\star_\Ms}\cdot_{f(\alpha_-)} a^{\star_{\alpha_-}})\otimes_\Zs (d^{\star_{\delta_-}}\cdot_{g(\delta_-)} y^{\star_\Ns} \cdot_{g(\gamma_-)} c^{\star_{\gamma_-}})
\\
&=b^{\star_{\beta_-}}\cdot_{f(\beta_-)}d^{\star_{\delta_-}}\cdot_{g(\delta_-)} 
(x\otimes_\Zs y)^{\star_\Ms\otimes_\Zs\star_\Ns}
\cdot_{f(\alpha_-)} a^{\star_{\alpha_-}}
\cdot_{g(\gamma_-)} c^{\star_{\gamma_-}}, 
\end{align*}
where we used the fact that $(\Ts_{\tau_\pm},\star_{\tau_\pm})=(\Ts_{(f,g)(\tau_\pm)},\star_{(f,g)(\tau_\pm)})$, for $\tau\in\{\alpha,\beta,\gamma,\delta \}$. 
\end{proof}

More generally, we can use the tensor product over subfamilies defined in remark~\ref{rem: otimes-tr}. 
\begin{remark}
Let $\left({}_{(\As_\alpha)_A}\Ms_{(\As_\beta)_B},((\star_\alpha)_A,\star_\Ms,(\star_\beta)_B)_{f_\Ms}\right)$ and $\left({}_{(\Bs_\gamma)_C}\Ns_{(\Bs_\delta)_D},((\star_\gamma)_C,\star_\Ns,(\star_\delta)_D)_{f_\Ns}\right)$ be two involutive multimodules. The ``internal tensor product'' $\Ms\otimes_\Gamma\Ns$ over an indexed family $\Gamma\subset \Sigma:=\{(\alpha,\beta) \ | \ \As_\alpha=\Bs_\beta \}$ of common subalgebras $\As_\beta=\Bs_{\gamma}$, with $(\beta,\gamma)\in \Gamma$, that is stable under $f:=(f_\Ms,f_\Ns)$, the disjoint union of the support involutions $f_\Ms,f_\Ns$:  
\begin{equation*}
(\xi,\zeta)\in\Gamma \imp (f(\xi),f(\zeta))\in\Gamma, \quad \forall \xi,\zeta\in A\uplus B\uplus C\uplus D, 
\end{equation*}
becomes an involutive multimodule with involution $\star:=\star_\Ms\otimes_\Gamma\star_\Ns$. The involution $\star$ is well-defined by universal factorization property of tensor products.\footnote{
Apart from checking directly that, under the stability condition, the involution is well-defined, it is also possible to obtain the same result, considering first the involution already defined in proposition~\ref{prop: inv-ot} and making use of  proposition~\ref{prop: inv-tr} together with remark~\ref{rem: inv-tr}.} 
\xqed{\lrcorner}
\end{remark}

\section{Pairing Dualities in $\Zs$-central Multimodules} \label{sec: multi-dual-paring}

Here we provide an extension, to the case of $\Zs$-central multimodules, of the notion of duality of vector spaces.  

\medskip 

Although tensor products are always introduced via their universal factorization property, and later used to provide examples of monoidal categories, in the literature duals are almost never defined via universal factorization properties and are rather described either with non-categorical definitions or as dual objects inside suitable monoidal categories. 

\medskip 

Our main purpose here will be to directly discuss the several pairing dualities for multi-modules.

\medskip

Let us more generally consider the case of $\Zs$-central left-$(\As_\alpha)_{\alpha\in A}$ right-$(\Bs_\beta)_{\beta\in B}$ multimodules ${}_{(\As_\alpha)_{\alpha\in A}}\Ms_{(\Bs_\beta)_{\beta\in B}}$. 
We can define several notions of duals, one for every subset of indexes $I\times J\subset A\times B$: 
\begin{definition}\label{def: dual-multi}
Given a $\Zs$-central multimodule  ${}_{(\As_\alpha)_{\alpha\in A}}\Ms_{(\Bs_\beta)_{\beta\in B}}$ over $\Rs_\Zs$-algebras $(\As_\alpha)_{\alpha\in A}$-$(\Bs_\beta)_{\beta\in B}$ and a family of indexes $I\times J\subset A\times B$, an \emph{($I,J$)-dual of the multimodule} $\Ms$ is a pair $(\Ns,\tau)$, where ${}_{(\Bs_\beta)_{\beta\in B}}\Ns_{(\As_\alpha)_{\alpha\in A}}$ is a \hbox{$\Zs$-central} $(\Bs_\beta)_{\beta\in B}$-$(\As_\alpha)_{\alpha\in A}$ multimodule over the $\Rs_\Zs$-algebras $(\Bs_\beta)_{\beta\in B}$-$(\As_\alpha)_{\alpha\in A}$ and 
$\tau:\Ns\times\Ms\to (\otimes^\Zs_{\alpha\in I}\As_\alpha)\otimes_\Zs(\otimes^\Zs_{\beta\in J}\Bs_\beta)$ is a $\Zs$-multilinear $(A-I,B-J)$-balanced $(I,J)$-multilinear map:\footnote{\label{foo: free-m}
With some abuse of notation, will denote by $\cdot$ the ``outer actions'' on the tensor product multimodule $(\otimes^\Zs_{\alpha\in I}\As_\alpha)\otimes_\Zs(\otimes^\Zs_{\beta\in J}\Bs_\beta)$ given by:  $a\cdot_{\beta_o}[(\otimes^\Zs_{\alpha\in I} x_\alpha)\otimes_\Zs(\otimes^\Zs_{\beta\in J} y_\beta)]\cdot_{\alpha_o} b
:=
(\otimes^\Zs_{\alpha\in I} x'_\alpha)\otimes_\Zs(\otimes^\Zs_{\beta\in J} y'_\beta)$, 
where $x'_\alpha:=\begin{cases}x_\alpha,\phantom{\ \ \cdot a} \quad \alpha\neq \alpha_o
\\ a\cdot x_{\alpha_o}, \quad \alpha={\alpha_o}
\end{cases}
$ and  $y'_\beta:=\begin{cases} 
y_\beta,\phantom{\ \ \cdot b} \quad \beta\neq \beta_o
\\ y_{\beta_o}\cdot b, \quad \beta={\beta_o}
\end{cases}
$  
for all $(\alpha_o,\beta_o), (\alpha,\beta)\in I\times J$, $(a,b)\in\As_{\alpha_o}\times\Bs_{\beta_o}$, $(x_\alpha,y_\beta)\in\As_\alpha\times\Bs_\beta$. 
} 
$\forall t,x\in\Ms$, 
\begin{gather*}
\tau(t,a\cdot_\alpha x\cdot_\beta b)=a\cdot_\alpha\tau(t,x)\cdot_\beta b, 
\quad 
\tau(b\cdot_\beta t\cdot_\alpha a,x)=b\cdot_\beta\tau(t,x)\cdot_\alpha a,  
\ \forall (a,b)\in\As_\alpha\times \Bs_\beta, \ (\alpha,\beta)\in I\times J, 
\\ 
\tau(b\cdot_\beta t\cdot_\alpha a,x)=\tau(t,a\cdot_\alpha x\cdot_\beta b), \quad \forall (a,b)\in \As_\alpha\times\Bs_\beta, \quad  (\alpha,\beta)\in(A-I)\times(B-J), 
\end{gather*}
satisfying the following universal factorization property: 
for any $(A-I,B-J)$-balanced $(I,J)$-multilinear map $\Phi:\widehat{\Ns}\times\Ms\to (\otimes^\Zs_{\alpha\in I}\As_\alpha)\otimes_\Zs(\otimes^\Zs_{\beta\in J}\Bs_\beta)$ 
where${}_{(\Bs_\beta)_{\beta\in B}}\widehat{\Ns}_{(\As_\alpha)_{\alpha\in A}}$ is another $\Zs$-central 
$(\Bs_\beta)_{\beta\in B}$-$(\As_\alpha)_{\alpha\in A}$ multimodule over $\Rs_\Zs$-algebras, 
there exists a unique morphism of multimodules $\hat{\Phi}:\widehat{\Ns}\to\Ns$ such that $\Phi=\tau\circ(\hat{\Phi},\id_\Ms)$. \footnote{
We are assuming here that, in the category of $\Zs$-central $\Rs_\Zs$-multimodules,   $\otimes^\Zs_{\alpha\in\varnothing}\As_\alpha:=\Rs=:\otimes^\Zs_{\beta\in\varnothing}\Bs_\beta$. 
} 
\end{definition}
Again, if an $(I,J)$-dual exists, it is unique up to a unique isomorphism of $(\Bs_\beta)_{\beta\in B}$-$(\As_\alpha)_{\alpha\in A}$ multimodules satisfying the previous universal factorization property. 
The existence is provided in the following result. 
\begin{theorem}
For every $(I,J)$ with $I\times J\subset A\times B$, there exists an $(I,J)$-dual $({}^{*_{I}}\Ms^{*_{J}},\tau)$ of the $\Zs$-central multimodule $\Ms$ over $\Rs_\Zs$-algebras $(\As_\alpha)_{\alpha\in A}$-$(\Bs_\beta)_{\beta\in B}$. 
\end{theorem}
\begin{proof}
For every $(I,J)$ with $I\times J\subset A\times B$, consider the following set: 
\begin{equation*}
{}^{*_{I}}\Ms^{*_{J}}:=\left\{\Ms\xrightarrow{\phi} (\otimes^\Zs_{\alpha\in I}\As_\alpha)\otimes_\Zs(\otimes^\Zs_{\beta\in J}\Bs_\beta) \ | \ \forall (\alpha,\beta)\in I\times J, \ \forall (a,b)\in\As_\alpha\times\Bs_\beta \st \phi(a\cdot_\alpha x\cdot_\beta b)=a\cdot_\alpha\phi(x)\cdot_\beta b \right\}. 
\end{equation*}
We see that ${}^{*_{I}}\Ms^{*_{J}}$ is a $\Zs$-central $(\Bs_\beta)_{\beta\in A}$-$(\As_\alpha)_{\alpha\in A}$ multimodule defining, 
for all $\phi,\psi\in {}^{*_{I}}\Ms^{*_{J}}$ and $x\in \Ms$:\footnote{
\label{foo: free-m2}
To avoid confusion, we will denote by $\bullet$ the ``inner actions'' on the tensor product multimodule $(\otimes^\Zs_{\alpha\in I}\As_\alpha)\otimes_\Zs(\otimes^\Zs_{\beta\in J}\Bs_\beta)$ given by:  $b\bullet_{\beta_o}[(\otimes^\Zs_{\alpha\in I} x_\alpha)\otimes_\Zs(\otimes^\Zs_{\beta\in J} y_\beta)]\bullet_{\alpha_o} a
:=
(\otimes^\Zs_{\alpha\in I} x'_\alpha)\otimes_\Zs(\otimes^\Zs_{\beta\in J} y'_\beta)$, 
where $x'_\alpha:=\begin{cases}x_\alpha,\phantom{\ \ \cdot a} \quad \alpha\neq \alpha_o
\\ x_{\alpha_o}\cdot a, \quad \alpha={\alpha_o}
\end{cases}
$ and  $y'_\beta:=\begin{cases} 
y_\beta,\phantom{\ \ \cdot b} \quad \beta\neq \beta_o
\\ b\cdot y_{\beta_o}, \quad \beta={\beta_o}
\end{cases}
$  
for all $(\alpha_o,\beta_o), (\alpha,\beta)\in I\times J$, $(a,b)\in\As_{\alpha_o}\times\Bs_{\beta_o}$, $(x_\alpha,y_\beta)\in\As_\alpha\times\Bs_\beta$. 
}
\begin{gather*}
\phi+\psi: x\mapsto \phi(x)+\psi(x), 
\\ 
b\bullet_\beta \phi\bullet_\alpha a: x\mapsto b\bullet_\beta \phi(x)\bullet_\alpha a, 
\quad \forall (\alpha,\beta)\in I\times J, \quad (a,b)\in \As_\alpha\times\Bs_\beta, 
\\
b\odot_\beta \phi\odot_\alpha a: x\mapsto \phi(a\cdot_\alpha x\cdot_\beta b),  
\quad \forall (\alpha,\beta)\in (A-I)\times (B-J), \quad (a,b)\in \As_\alpha\times\Bs_\beta. 
\end{gather*}
The evaluation map $\tau(\phi,x):=\phi(x)$, for all $\phi\in {}^{*_{I}}\Ms^{*_{J}}$ and $x\in \Ms$ turns out to be an $(A-I,B-J)$-balanced $(I,J)$-multilinear map $\tau:{}^{*_{I}}\Ms^{*_{J}}\times\Ms\to 
(\otimes^\Zs_{\alpha\in I}\As_\alpha)\otimes_\Zs(\otimes^\Zs_{\beta\in J}\Bs_\beta)$. 

\medskip 

To every $\Zs$-multilinear $(A-I)$-$(B-J)$-balanced and $(I,J)$-multilinear map 
$\Phi:\widehat{\Ns}\times\Ms \to (\otimes^\Zs_{\alpha\in I}\As_\alpha)\otimes_\Zs(\otimes^\Zs_{\beta\in J}\Bs_\beta)$, the usual Curry isomorphism associates the map $\widehat{\Phi}:\widehat{\Ns} \to 
[(\otimes^\Zs_{\alpha\in I}\As_\alpha)\otimes_\Zs(\otimes^\Zs_{\beta\in J}\Bs_\beta)]^\Ms$ that to every element $t\in\widehat{\Ns}$ associates the map 
$\widehat{\Phi}_t:\Ms\to (\otimes^\Zs_{\alpha\in I}\As_\alpha)\otimes_\Zs(\otimes^\Zs_{\beta\in J}\Bs_\beta)$ given by $\widehat{\Phi}_t(x):=\Phi(t,x)$, for all $x\in\Ms$. The defining properties of $\Phi$ entail that $\widehat{\Phi}_t\in{}^{*_I}\Ms^{*_J}$, for all $t\in\widehat{\Ns}$ and that the map $\widehat{\Phi}:\widehat{\Ns}\to{}^{*_I}\Ms^{*_J}$ given by $t\mapsto \widehat{\Phi}_t$ is a morphism of $(\Bs_\beta)_{\beta\in B}$-$(\As_\alpha)_{\alpha\in A}$ multimodules. Finally $\Phi(t,x)=\widehat{\Phi}_t(x)=\tau(\widehat{\Phi}_t,x)$, for $t\in\widehat{\Ns}$ and $x\in\Ms$. 
\end{proof}

For every pair of families of unital associative $\Zs$-central $\Rs_\Zs$-algebras $(\As_\alpha)_{\alpha\in A}$ and $(\Bs_\beta)_{\beta\in B}$ and every family of indexes $I\times J\subset A\times B$,  $(I,J)$-transposition functors (and evaluation natural transformations) give us a \textit{contravariant right semi-adjunction} according to the definitions fully recalled in appendix~\ref{sec: functorial-pairing}, remark~\ref{rem: s-aj}. 
\begin{theorem}\label{th: sadj}
Let ${}_{(\As_\alpha)_A}\Mf_{(\Bs_\beta)_B}$ be the category with objects $\Zs$-central $(\As_\alpha)_A$-$(\Bs_\beta)_B$ multimodules over unital associative \hbox{$\Rs$-algebras} and with morphism \hbox{$\Zs$-lin}\-ear maps $\Ms_1\xrightarrow{\Phi}\Ms_2$ such that $\Phi(a\cdot_\alpha x\cdot_\beta b)=a\cdot_\alpha \Phi(x)\cdot_\beta b$, for all $(\alpha,\beta)\in A\times B$, $(a,b)\in\As_\alpha\times\Bs_\beta$, $x\in\Ms$.
 
\medskip 

For every subset $I\times J\subset A\times B$ of indexes, we have a different contravariant right semi-adjoint functorial pairing $\quad \underline{{}_I\flat_J} \ | \stackrel[\theta]{\vartheta}{\leftrightarrows}|\ {}_I\sharp_J \quad$ between the \emph{transposition functors}\footnote{
The apparent distinction between $\flat$ and $\sharp$ is purely formal since they interchange by permuting the sets of indexes: ${}_{I}\flat_{J}={}_{J}\sharp_{I}$. 
}  
$\xymatrix{{}_{(\As_\alpha)_A}\Mf_{(\Bs_\beta)_B} \rtwocell^{{}_{{}_I\flat_J}}_{{}^{{}_I\sharp_J}}{'} & {}_{(\Bs_\beta)_B}\Mf_{(\As_\alpha)_A}}$ that 
\begin{itemize}
\item
on objects of the respective categories, are given by duals:  
\begin{equation*}
{}_{I}\flat_{J}:\Ms\mapsto {}^{*_{I}}\Ms^{*_{J}}, \quad \forall \Ms\in \Ob_{{}_{(\As_\alpha)_A}\Mf_{(\Bs_\beta)_B}}, 
\quad \quad 
{}_{I}\sharp_{J}:\Ns\mapsto {}^{*_{J}}\Ns^{*_{I}}, \quad \forall \Ns\in \Ob_{{}_{(\Bs_\beta)_B}\Mf_{(\As_\alpha)_A}};  
\end{equation*}
\item
on morphisms $(\Ms_2\xleftarrow{\mu}\Ms_1)\in\Hom_{{}_{(\As_\alpha)_A}\Mf_{(\Bs_\beta)_B}}(\Ms_1;\Ms_2)$ and   $(\Ns_2\xleftarrow{\nu}\Ns_1)\in\Hom_{{}_{(\Bs_\beta)_B}\Mf_{(\As_\alpha)_A}}(\Ns_1;\Ns_2)$, are respectively given by $\mu$-pull-backs and $\nu$-pull-backs: 
\begin{gather} \label{eq: transp}
(\Ms_2)^{{}_{I}\flat_{J}}\xrightarrow{\mu^{{}_{I}\flat_{J}}}(\Ms_1)^{{}_{I}\flat_{J}}, 
\quad \mu^{{}_{I}\flat_{J}}(\phi):=\phi\circ \mu, \quad \forall \phi\in {}^{*_I}\Ms_2^{*_J}, 
\\ \notag
(\Ns_2)^{{}_{I}\sharp_{J}}\xrightarrow{\nu^{{}_{I}\sharp_{J}}}(\Ns_1)^{{}_{I}\sharp_{J}}, 
\quad \nu^{{}_{J}\flat_{I}}(\psi):=\psi\circ \nu,
\quad \forall \psi\in {}^{*_J}\Ns_2^{*_I};  
\end{gather}
\end{itemize}
where unit and co-unit of the semi-adjunction are given by the following \emph{natural evaluation transformations}:\footnote{
These evaluations maps are just obtained applying Curry isomorphism to the pairing duality $\tau$ in definition~\ref{def: dual-multi}.
}
\begin{gather*} 
\id_{{}_{(\As_\alpha)_{\alpha\in A}}\Mf_{(\Bs_\beta)_{\beta\in B}}}\xrightarrow{\theta}\sharp\circ\flat, 
\quad \Ms\mapsto \theta^\Ms, \quad 
\Ms \xrightarrow{\theta^\Ms} (\Ms^{{}_{I}\flat_{J}})^{{}_{J}\sharp_{I}}, 
\quad \theta^\Ms_x:\phi\mapsto \phi(x), \quad \forall \phi\in \Ms^{{}_{I}\flat_{J}},\ x\in \Ms, 
\\
\id_{{}_{(\Bs_\beta)_{\beta\in B}}\Mf_{(\As_\alpha)_{\alpha\in A}}}\xrightarrow{\vartheta}\flat\circ\sharp, 
\quad 
\Ns\mapsto \vartheta^\Ns, 
\quad 
\Ns\xrightarrow{\vartheta^\Ns} (\Ns^{{}_{J}\sharp_{I}})^{{}_{I}\flat_{J}}, 
\quad \vartheta^\Ns_y:\psi\mapsto \psi(y), \quad \forall \psi\in \Ns^{{}_{J}\sharp_{I}},\ y\in \Ns. 
\end{gather*}
Restricting the previous contravariant right semi-adjunction $\ \underline{{}_I\flat_J} \ | \stackrel[\theta]{\vartheta}{\leftrightarrows}|\ {}_I\sharp_J \ $ to the full reflective subcategories (whose objects are those multimodules for which the evaluation natural transformations are isomorphisms), we obtain a categorical duality. 
\end{theorem}
\begin{proof}
The contravariant functorial nature of ${}_I\flat_J$ and ${}_I\sharp_J$ is standard from their definitions.  

\medskip 

By direct computation $\theta^\Ms$ is a morphism in ${}_{(\As_\alpha)_{\alpha\in A}}\Mf_{(\Bs_\beta)_{\beta\in B}}$ and $\vartheta^\Ns$ is a morphism in ${}_{(\Bs_\beta)_{\beta\in B}}\Mf_{(\As_\alpha)_{\alpha\in A}}$ furthermore for every pair of morphisms $(\Ms_1)\xrightarrow{\mu}(\Ms_2)$ in ${}_{(\As_\alpha)_{\alpha\in A}}\Mf_{(\Bs_\beta)_{\beta\in B}}$ and  $(\Ns_1)\xrightarrow{\nu}(\Ns_2)$ in ${}_{(\Bs_\beta)_{\beta\in B}}\Mf_{(\As_\alpha)_{\alpha\in A}}$: 
\begin{equation*}
\theta^{\Ms_2}\circ \mu=(\mu^{{}_{I}\flat_{J}})^{{}_{J}\sharp_{I}} \circ \theta^{\Ms_1}, 
\quad \quad 
\vartheta^{\Ns_2}\circ \nu=(\nu^{{}_{J}\sharp_{I}})^{{}_{I}\flat_{J}} \circ \vartheta^{\Ns_1}. 
\end{equation*}
Finally we check the right semi-adjunction condition $\underline{{}_I\flat_J} \ | \stackrel[
]{
}{\leftrightarrows}|\ {}_I\sharp_J$ using formula~\eqref{eq: r-s-aj}: 
\begin{align*}
({{}_I\flat_J}(\theta^\Ms))\circ \vartheta^{(\Ms^{{}_{I}\flat_{J}})}&=\iota_{(\Ms^{{}_{I}\flat_{J}})}, 
\quad \forall \Ms\in\Ob_{{}_{(\As_\alpha)_A}\Mf_{(\Bs_\beta)_B}}, 
\\ 
[({{}_I\flat_J}(\theta^\Ms)\circ \vartheta^{\Ms^{{}_{I}\flat_{J}}})(\phi)](x)
&
=[(\theta^\Ms)^{{}_I\flat_J}(\vartheta^{\Ms^{{}_{I}\flat_{J}}}_\phi)](x)
=[(\vartheta^{\Ms^{{}_{I}\flat_{J}}}_\phi)\circ \theta^\Ms](x)
= \vartheta^{\Ms^{{}_{I}\flat_{J}}}_\phi(\theta^\Ms_x)
\\
&
=\theta_x^\Ms(\phi)
=\phi(x)
=[\iota_{(\Ms^{{}_{I}\flat_{J}})}(\phi)](x), 
\quad \quad \forall \phi\in \Ms^{{}_{I}\flat_{J}}, \quad x\in\Ms. 
\end{align*} 
For the full reflective subcategories of the semi-adjunction we have a categorical duality (see remark~\ref{rem: s-aj}). 
\end{proof}

\begin{remark}
There is of course the possibility to define also \emph{${}_I\gamma_J$-conjugate duals} of ${}_{(\As_\alpha)_{A}}\Ms_{(\Bs_\beta)_B}$ for any family of $\Rs_\Zs$-conjugations $(\gamma_k)_{k\in I\uplus J}$, for 
$(i,j)\in I\times J\subset A\times B$. For this purpose is just enough to repeat the previous construction of duals utilizing maps that are $\gamma_k$-conjugate-$\Rs_\Zs$-linear. Whenever $\gamma_k=\id_{\Rs}$, for all $k\in I\uplus J$, we re-obtain the previous definition as a special case. 
\xqed{\lrcorner}
\end{remark}

Here we discuss how our definition of duals relates to already available notions in the case of bimodules. 
\begin{remark} 
The notion of $I$-$J$ dual of an $(\As_\alpha)_A$-$(\Bs_\beta)_B$ multimodule over $\Rs_\Zs$-algebras that we have just introduced in definition~\ref{def: dual-multi} is a direct generalization of some much more familiar constructs for bimodules.  

\medskip 

Here below, we consider an $\As$-bimodule ${}_\As\Ms_\As$ as an ${(\As_\alpha)_A}$-$(\Bs_\beta)_B$-multimodule, with 
$A:=\{\alpha_o\}$, $B:=\{\beta_o \}$, $A\times B=\{(\alpha_o,\beta_o)\}$ singleton sets and with $\As_{\alpha_o}:=\As=:\Bs_{\beta_o}$. 

\medskip 

The ``double dual''\footnote{
To be precise, the double dual is obtained choosing here $\Zs:=\KK$; this is a slight generalization that we found particularly useful in our treatment of contravariant non-commutative differential calculus (see footnote~\ref{foo: ncdc}).
} 
$\Ms^\vee$ of an \hbox{$\As$-bimodule} ${}_\As\Ms_\As$ (see for example~\cite[section~2.1]{Fer17} for more details) is the central $\Zs$-bimodule $\Ms^\vee:=\Hom_{{}_\As\Mf_\As}(\Ms;{}_{\cdot \ }(\As\otimes_\Zs\As)_{\ \cdot})$ of covariant homomorphisms of bimodules, from ${}_\As\Ms_\As$, with values into ${}_{\cdot \ }(\As\otimes_\Zs\As)_{\ \cdot}$ seen as an $\As$-bimodule with the ``exterior actions'' given by: $a\cdot (x\otimes_\Zs y)\cdot b:=(ax)\otimes_\Zs(yb)$, for all $x,y,a,b\in \As$; where $\Ms^\vee$ is an $\As$-bimodule with the actions $(b\odot \phi\odot a)(x):=\phi(a\cdot x\cdot b)$, for all $a,b\in\As$, $x\in \Ms$, $\phi\in\Ms^\vee$. 
Taking $I\times J = A\times B$ as a singleton (only one right and only one left action) in definition~\ref{def: dual-multi},  we see that $\Ms^\vee={}^{*_I}\Ms^{*_J}$. 

\medskip 

The well-known notions (see for example~\cite{Bo97}) of ``right dual'' $\Ms^*:=\Hom_\As(\Ms_\As;\As)$ and ``left dual'' ${}^*\Ms:=\Hom_\As({}_\As\Ms;\As)$ of a bimodule ${}_\As\Ms_\As$ are just the dual of the right $\As$-module $\Ms_\As$ (respectively the dual of the left $\As$-module ${}_\As\Ms$), as in~\cite[section~II.3]{Bou89}, equipped with the following actions $(a\cdot \phi\odot b)(x):=a\phi(bx)$, for all $a,b\in\As$, $x\in\Ms$ and $\phi\in\Ms^*$ (respectively $(a\odot \phi\cdot b)(x):=\psi(xa)b$, for all $a,b\in\As$, $x\in\Ms$ and $\psi\in{}^*\Ms$). 
When $A\times B$ is a singleton, taking $I:=\varnothing$, $J:=B$, we recover $\Ms^*={}^{*_I}\Ms^{*_J}$ and, when $I:=B$, $J:=\varnothing$, we get ${}^*\Ms={}^{*_I}\Ms^{*_J}$. 

\medskip 

Finally the ``scalar dual'' of a bimodule\footnote{
For algebras over $\Rs:=\Zs:=\KK$ this is just the usual dual as a $\KK$-vector space. 
} ${}_\As\Ms_\As$, defined as $\Ms':=\Hom_\Rs(\Ms;\Rs\otimes_\Zs\Rs)$, equipped with the actions $(b\odot \phi\odot a)(x):=\phi(a\cdot x\cdot b)$, for all $a,b\in\As$, $x\in\Ms$ and $\phi\in\Ms'$, can be obtained from our definition as $\Ms'={}^{*_I}\Ms^{*_J}$, taking $A\times B$ to be, as usual, a singleton and $I:=\varnothing=:J$.  
\xqed{\lrcorner}
\end{remark}

In the following we study the ``inclusion relations'' between the different duals of a given multimodule. 
\begin{remark}
Consider the auxiliary $(I_2,J_2)$-global $(I_1,J_1)$-dual multimodules, for $I_1\times J_1\subset I_2\times J_2\subset A\times B$: 
\begin{equation*}
{}^{I_1}_{I_2}\Ms^{J_1}_{J_2}:=
\left\{\Ms\xrightarrow{\phi} (\otimes^\Zs_{\alpha\in I_2}\As_\alpha)\otimes_\Zs(\otimes^\Zs_{\beta\in J_2}\Bs_\beta) \ | \ \forall (\alpha,\beta)\in I_1\times J_1, \ \forall (a_\alpha,b_\beta)\in\As_\alpha\times\Bs_\beta \st \phi(a_\alpha xb_\beta)=a_\alpha\phi(x)b_\beta \right\}, 
\end{equation*}
equipped with the multimodule actions specified as follows (see footnotes~\ref{foo: free-m} \ref{foo: free-m2}), for all $\phi\in{}^{I_1}_{I_2}\Ms^{J_1}_{J_2}$ and $x\in\Ms$: 
\begin{gather}\notag 
b\bullet_\beta \phi\bullet_\alpha a: x\mapsto b\bullet_\beta \phi(x)\bullet_\alpha a,  
\quad \forall (\alpha,\beta)\in I_2\times J_2, \quad (a,b)\in \As_\alpha\times\Bs_\beta, 
\\ \notag
b\odot_\beta \phi\odot_\alpha a: x\mapsto \phi(a\cdot_\alpha x\cdot_\beta b), 
\quad \forall (\alpha,\beta)\in (A-I_1)\times (B-J_1), \quad (a,b)\in \As_\alpha\times\Bs_\beta,  
\\ \label{eq: extra}
(a\cdot_\alpha \phi\cdot_\beta b):=x \mapsto a\cdot_\alpha\phi(x)\cdot_\beta b, \quad \forall (\alpha,\beta)\in (I_2-I_1)\times(J_2-J_1). 
\end{gather}
Notice that whenever $(I_1,J_1)=(I_2,J_2)$, we have ${}^{I_1}_{I_2}\Ms^{J_1}_{J_2}={}^{*_{I_1}}\Ms^{*_{J_1}}$ as a multimodule and that the extra multimodule actions in line~\eqref{eq: extra} appear only when $I_1\times J_1\neq I_2\times J_2$.  

\medskip 

If $I_1\times J_1\subset I_1'\times J_1'$ we have natural set theoretic inclusions: $\xymatrix{{}^{I'_1}_{I_2}\Ms^{J'_1}_{J_2} \ar@{^(->}[r]^{\eta_\Ms} & {}^{I_1}_{I_2}\Ms^{J_1}_{J_2}}$, that are also covariant morphisms of multimodules for all the common actions involved (inner target actions $\bullet$ for indexes in $I_2\times J_2$; internal source actions $\odot$ for indexes in $(A-I_1')\times(B-J_1')$ and external target action 
$\cdot \ $ for indexes in the set $(I_2-I_1')\times(J_2-J_1')$). 
Keeping $(I_1,J_1)$ fixed, if $I_2\times J_2\subset I_2'\times J_2'$, we define the following embedding map: 
\begin{equation*}
{}^{I_1}_{I_2}\Ms^{J_1}_{J_2}\xrightarrow{\quad \zeta_\Ms:=(\otimes^\Zs_{\alpha\in I_2'-I_2}1_{\As_\alpha})\otimes_\Zs - \otimes_\Zs(\otimes^\Zs_{\beta\in J_2'-J_2}1_{\Bs_\beta})\quad}{}^{I_1}_{I'_2}\Ms^{J_1}_{J'_2}, 
\end{equation*}
that to every $\phi\in{}^{I_1}_{I_2}\Ms^{J_1}_{J_2}$ associates the map $\Ms\ni x\mapsto (\otimes^\Zs_{\alpha\in I_2'-I_2}1_{\As_\alpha})\otimes_\Zs \phi(x) \otimes_\Zs(\otimes^\Zs_{\beta\in J_2'-J_2}1_{\Bs_\beta})$ in ${}^{I_1}_{I'_2}\Ms^{J_1}_{J'_2}$, that is actually a covariant morphism of multimodules for all the common relevant actions involved (inner target actions $\bullet$ for indexes in $I_2\times J_2$; internal source actions $\odot$ for indexes in $(A-I_1)\times(B-J_1)$ and external target action $\cdot \ $ for indexes in $(I_2-I_1)\times(J_2-J_1)$).  
\xqed{\lrcorner}
\end{remark}

\begin{proposition}
Given a multimodule ${}_{(\As_\alpha)_{\alpha\in A}}\Ms_{(\Bs_\beta)_{\beta\in B}}$, for all the inclusions $I_1\times J_1\subset I_2\times J_2\subset A\times B$ of indexes, we have the following natural transformations between contravariant functors from the category ${}_{(\Bs_\beta)_{\beta\in B}}\Mf_{(\As_\alpha)_{\alpha\in A}}$ into the category $\Mf_{[\Rs_\Zs]}$ of $\Zs$-central multimodules over $\Rs_\Zs$-algebras:  
\begin{equation*}
\Ms \quad \mapsto \quad  \left[ 
{}^{*_{I_2}}\Ms^{*_{J_2}} 
\xrightarrow{\tiny 
\xymatrix{{}^{I'_1}_{I_2}\Ms^{J'_1}_{J_2} \ar@{^(->}[r]^{\eta_\Ms}& {}^{I_1}_{I_2}\Ms^{J_1}_{J_2}}
}
{}^{I_1}_{I_2}\Ms^{J_1}_{J_2} 
\xleftarrow{\zeta_\Ms:=(\otimes^\Zs_{\beta\in J-J'}1_{\Bs_\beta})\otimes_\Zs -\otimes_\Zs(\otimes^\Zs_{\alpha\in I-I'}1_{\As_\alpha})} {}^{*_{I_1}}\Ms^{*_{J_1}}
\right]. 
\end{equation*}
\end{proposition}
\begin{proof}
The passage associating to a multimodule $\Ms\in\Ob_{{}_{(\As_\alpha)_A}\Mf_{(\Bs_\beta)_B}}$ its $(I_2,J_2)$-global $(I_1,J_1))$-dual multimodule ${}^{I_1}_{I_2}\Ms^{J_1}_{J_2}$, 
is a contravariant functor acting on morphisms by transposition as in equation~\eqref{eq: transp}: 
\begin{equation*}
\left(\Ms\xrightarrow{\mu}\Ns \right) \quad \mapsto \left( {}^{I_1}_{I_2}\Ms^{J_1}_{J_2} \xleftarrow{\mu^\bullet} {}^{I_1}_{I_2}\Ns^{J_1}_{J_2} \right), \quad \text{where} \quad \mu^\bullet(\phi):=\phi\circ \mu, 
\end{equation*}
and from $\eta_\Ms\circ\mu^{{}_{I_2}\flat_{J_2}}=\mu^\bullet\circ \eta_\Ns$ and $\mu^\bullet\circ\zeta_\Ns=\zeta_\Ms\circ\mu^{{}_{I_1}\flat_{J_1}}$ we see that $\eta$ and $\zeta$ are natural transformations. 
\end{proof}

\begin{theorem}
For every inclusion of indexes $I_1\times J_1\subset I_2\times J_2\subset A\times B$, considering the two contravariant right semi-adjunctions 
$\left[\underline{{}_{I_k}\flat_{J_k}} \ | \stackrel[{}^{I_k}\theta^{J_k}]{{}^{I_k}\vartheta^{J_k}}{\leftrightarrows}|\ {}_{I_k}\sharp_{J_k} \right]$, for $k=1,2$, as in theorem~\ref{th: sadj}, we have a \emph{morphism of contravariant right semi-adjunctions} defined in the following way: 
\begin{itemize}
\item 
for all morphisms $\Ms\xrightarrow{\mu}\Ns$ in ${}_{(\As_\alpha)_A}\Mf_{(\Bs_\beta)_B}$ 
and $\Ps\xrightarrow{\nu}\Qs$ in ${}_{(\Bs_\beta)_B}\Mf_{(\As_\alpha)_A}$ 
we have commutative diagrams: 
\begin{equation*}
\xymatrix{
 & (\Ps)^{{}_{I_1}\sharp_{J_1}} \ar[dl]_{\zeta_\Ps} & \ar[l]_{\quad \nu^{{}_{I_1}\sharp_{J_1}} } (\Qs)^{{}_{I_1}\sharp_{J_1}}\ar[dr]^{\zeta_\Qs} &  
\\
{}^{I_1}_{I_2}\Ps^{J_1}_{J_2} & & & {}^{I_1}_{I_2}\Qs^{J_1}_{J_2} \ar[lll]_{\nu^\bullet}
\\
 & (\Ps)^{{}_{I_2}\sharp_{J_2}} \ar[ul]^{\eta_\Ps} & \ar[l]_{\quad \nu^{{}_{I_2}\sharp_{J_2}} } (\Qs)^{{}_{I_2}\sharp_{J_2}} \ar[ur]_{\eta_\Qs} &  
}
\quad \quad 
\xymatrix{
 & (\Ms)^{{}_{I_1}\flat_{J_1}} \ar[ld]_{\zeta_\Ms} & \ar[l]_{\quad \mu^{{}_{I_1}\flat_{J_1}} } (\Ns)^{{}_{I_1}\flat_{J_1}} \ar[dr]^{\zeta_\Ns} &  
\\
{}^{I_1}_{I_2}\Ms^{J_1}_{J_2} & & & {}^{I_1}_{I_2}\Ns^{J_1}_{J_2} \ar[lll]_{\mu^\bullet}
\\
 & (\Ms)^{{}_{I_2}\flat_{J_2}} \ar[ul]^{\eta_\Ms}  & \ar[l]_{\quad \mu^{{}_{I_2}\flat_{J_2}} } (\Ns)^{{}_{I_2}\flat_{J_2}} \ar[ur]_{\eta_\Ns}  &  
}
\end{equation*}
\item 
for every pair of objects $\Ms$ in ${}_{(\As_\alpha)_A}\Mf_{(\Bs_\beta)_B}$ and $\Qs$ in  ${}_{(\Bs_\beta)_B}\Mf_{(\As_\alpha)_A}$ we have the commuting diagrams: 
\begin{equation*}
\xymatrix{
 & ((\Ms)^{{}_{I_1}\flat_{J_1}})^{{}_{I_1}\sharp_{J_1}} \ar[rr]^{\zeta_{(\Ms)^{{}_{I_1}\flat_{J_1}}}}
&  & {}^{J_1}_{J_2}((\Ms)^{{}_{I_1}\flat_{J_1}})^{I_1}_{I_2}
\\
\Ms \ar[ur]^{{{}^{I_1}\theta^{J_1}}_\Ms} \ar[dr]_{{{}^{I_2}\theta^{J_2}}_\Ms} \ar[rrr]_{\ev^\Ms \quad \quad} & & & {}^{J_1}_{J_2}({}^{I_1}_{I_2}\Ms^{J_1}_{J_2})^{I_1}_{I_2} \ar[u]_{\zeta_\Ms^\bullet} \ar[d]^{\eta_\Ms^\bullet} \ ,
\\
 & ((\Ms)^{{}_{I_2}\flat_{J_2}})^{{}_{I_2}\sharp_{J_2}}   \ar[rr]_{\eta_{(\Ms)^{{}_{I_1}\flat_{J_2}}}}
&  & {}^{J_1}_{J_2}((\Ms)^{{}_{I_2}\flat_{J_2}})^{I_1}_{I_2}
}
\quad \quad 
\xymatrix{
 & ((\Qs)^{{}_{I_1}\sharp_{J_1}})^{{}_{I_1}\flat_{J_1}} \ar[rr]^{\zeta_{(\Qs)^{{}_{I_1}\sharp_{J_1}}}}
&  & {}^{J_1}_{J_2}((\Qs)^{{}_{I_1}\sharp_{J_1}})^{I_1}_{I_2}
\\
\Qs \ar[ur]^{{{}^{I_1}\vartheta^{J_1}}_\Qs} \ar[dr]_{{{}^{I_2}\vartheta^{J_2}}_\Qs}\ar[rrr]_{\ev^\Qs\quad \quad } & & & {}^{J_1}_{J_2}({}^{I_1}_{I_2}\Qs^{J_1}_{J_2} )^{I_1}_{I_2}  \ar[u]_{\zeta_\Qs^\bullet} \ar[d]^{\eta_\Qs^\bullet}.
\\
 & ((\Qs)^{{}_{I_2}\sharp_{J_2}})^{{}_{I_2}\flat_{J_2}}   \ar[rr]_{\eta_{(\Qs)^{{}_{I_2}\sharp_{J_2}}}}
&  & {}^{J_1}_{J_2}((\Qs)^{{}_{I_2}\sharp_{J_2}})^{I_1}_{I_2}
}
\end{equation*}
\end{itemize}

\medskip 

Considering the category $\If$ (actually the poset) of pairs $(I,J)$, with $I\times J\subset A\times B$, where for every inclusion $I_1\times J_1\subset I_2\times J_2$ there is a unique morphism of pairs $(I_1,J_1)\xrightarrow{} (I_2,J_2)$, 
we have that every multimodule $\Ms$ has an associated \emph{dual functor} $\If\xrightarrow{\bigstar_\Ms}\Sf$ into contravariant right semi-adjunctions: $\bigstar_\Ms:(I,J)\mapsto \left[\underline{{}_I\flat_J} \ | \stackrel[{}^I\theta^J]{{}^I\vartheta^J}{\leftrightarrows}|\ {}_I\sharp_J \right]$. 
\end{theorem}
\begin{proof} 
In the first pair of diagram, due to the exchange symmetry $\mu\leftrightarrow \nu$ it is sufficient to prove the second. 
Taking $\phi\in(\Ns)^{{}_{I_1}\flat_{J_1}}$ and $\psi\in(\Ns)^{{}_{I_2}\flat_{J_2}}$ we immediately get, for $x\in\Ms$: 
\begin{align*}
[\zeta_\Ms\circ \mu^{{}_{I_1}\flat_{J_1}} (\phi)](x)
&= [\zeta_\Ms(\phi\circ\mu)](x)
=(\otimes^\Zs_{\alpha\in I_2-I_1}1_{\As_\alpha})
\otimes_\Zs\phi(\mu(x))\otimes_\Zs(\otimes^\Zs_{\beta\in J_2-J_1}1_{\Bs_\beta})
\\
&=[\zeta_\Ns (\phi)\circ \mu](x)
=[\mu^\bullet\circ\zeta_\Ns (\phi)](x), 
\\
[\eta_\Ms\circ \mu^{{}_{I_2}\flat_{J_2}} (\psi)](x)
&= [\eta_\Ms(\psi\circ\mu)](x)
=[\eta_\Ns (\psi)\circ \mu](x)
=\phi(\mu(x))
=[\mu^\bullet\circ\eta_\Ns (\psi)](x).
\end{align*}

\medskip 

In the second pair of commuting diagrams, by the exchange symmetry $\Ms\leftrightarrow \Qs$, it is enough to prove the first.
Consider $x\in\Ms$, $\phi\in (\Ms)^{{}_{I_1}\flat_{J_1}}$ and $\psi\in(\Ms)^{{}_{I_2}\flat_{J_2}}$: 
\begin{align*}
[\zeta_\Ms^\bullet\circ\ev^\Ms(x)](\phi) 
&=[\zeta_\Ms^\bullet(\ev^\Ms_x)](\phi)
=\ev^\Ms_x(\zeta_\Ms(\phi))
=(\otimes_{\alpha\in I_2-I_1}^\Zs 1_{\As_\alpha})\otimes_\Zs\phi(x)\otimes_\Zs(\otimes_{\beta\in J_2-J_1}^\Zs 1_{\Bs_\beta})
\\
&=[\zeta_{(\Ms)^{{}_{I_1}\flat_{J_1}}} (({}^{I_2}\theta^{J_2}_\Ms)_x)](\phi), 
\\
[\eta_\Ms^\bullet\circ\ev^\Ms(x)](\psi) 
&=[\eta_\Ms^\bullet(\ev^\Ms_x)](\psi)
=\ev^\Ms_x(\eta_\Ms(\psi))
=\psi(x)
=[\eta_{(\Ms)^{{}_{I_1}\flat_{J_2}}}({}^{I_2}\theta^{J_2}_\Ms)_x](\psi).
\end{align*}

We define a poset category $\If$ of index pairs via the order relation $(I_1,J_1)\leq(I_2,J_2) \iff (I_1\subset I_2)\wedge (J_1\subset J_2)$. 

\medskip 

We consider $\Sf_\Ms$ the category whose objects are contravariant right semi-adjunctions and whose morphisms are specified by the previous commuting diagrams of natural transformations $\zeta,\eta$. 

\medskip 

To every index pair $(I,J)\in\Ob_\If$ we associate the contravariant right semi-adjunction $\bigstar_\Ms^{(I,J)}:=\left[\underline{{}_I\flat_J} \ | \stackrel[{}^I\theta^J]{{}^I\vartheta^J}{\leftrightarrows}|\ {}_I\sharp_J \right]$
and to every morphism $(I_1,J_1)\leq(I_2,J_2)$ in $\If$ we associate the morphism $\bigstar_\Ms^{(I_2,J_2)}\xrightarrow{(\zeta,\eta)}\bigstar_\Ms^{(I_1,J_1)}$ of contravariant right semi-adjunctions. 
We notice that $\If\xrightarrow{\bigstar_\Ms}\Sf_\Ms$ is a contravariant functor. 
\end{proof}

\section{Traces and Inner Products on Multimodules} \label{sec: traces-ip}

In the first part of this section we generalize to the setting of multimodules the well-known multilinear algebraic operations producing contractions of tensors (over pairs of contravariant/covariant indexes) over a vector space and hence the equally familiar notion of trace of linear operators.

\medskip 

We proceed again introducing the relevant universal factorization properties. 

\begin{definition} 
Given a $\Zs$-central multimodule ${}_{(\As_\alpha)_A}\Ms_{(\As_\beta)_B}$ over $\Rs_\Zs$-algebras and $\Gamma\subset (A\uplus B)\times (A\uplus B)$ an injective symmetric relation\footnote{We can also assume that $\Gamma$ is irreflexive: $(\xi,\zeta)\in\Gamma \imp \xi\neq \zeta$; since ``tracing an action over itself'' does not have any effect.} such that $\As_\xi=\As_\zeta$, for all $(\xi,\zeta)\in\Gamma$, let  $A^\Gamma:=A-\dom(\Gamma)$ and $B^\Gamma:=B-\im(\Gamma)$. 

\medskip 

A $\Zs$-linear map ${}_{(\As_\alpha)_A}\Ms_{(\As_\beta)_B}\xrightarrow{T}\Vs$, of $\Zs$-central bimodules, is \emph{$\Gamma$-tracial} if it satisfies the following properties: 
\begin{align*}
T(a\cdot_\xi x)&=T(a\cdot_\zeta x), \quad \forall x\in\Ms, \ \forall a\in \As_\xi=\As_\zeta, \quad 
(\xi,\zeta)\in\Gamma\cap (A\times A), 
\\
T(x\cdot_\xi a)&=T(x\cdot_\zeta a), \quad \forall x\in\Ms, \ \forall a\in \As_\xi=\As_\zeta, \quad 
(\xi,\zeta)\in\Gamma\cap (B\times B),  
\\
T(a\cdot_\xi x)&=T(x\cdot_\zeta a), \quad \forall x\in\Ms, \ \forall a\in \As_\xi=\As_\zeta, \quad 
(\xi,\zeta)\in\Gamma\cap (A\times B).
\end{align*}

\medskip 

A \emph{$\Gamma$-contraction of the multimodule} ${}_{(\As_\alpha)_A}\Ms_{(\As_\beta)_B}$, consists of a $\Gamma$-tracial 
morphism of \hbox{$(\As_\alpha)_{\alpha\in A^\Gamma}$-$(\As_\beta)_{\beta\in B^\Gamma}$} multimodules
over $\Rs_\Zs$-algebras 
${}_{(\As_\alpha)_{A^\Gamma}}\Ms_{(\As_\beta)_{B^\Gamma}}\xrightarrow{T^\Gamma_\Ms} {}_{(\As_\alpha)_{A^\Gamma}}\Ms\ |\ \Gamma_{(\As_\beta)_{B^\Gamma}}$
such that the following universal factorization property is satisfied: 
for any $\Gamma$-tracial morphism 
${}_{(\As_\alpha)_{A^\Gamma}}\Ms_{(\As_\beta)_{B^\Gamma}}\xrightarrow{T} {}_{(\As_\alpha)_{A^\Gamma}}\Ns_{(\As_\beta)_{B^\Gamma}}$
of multimodules, there exists a unique morphism 
${}_{(\As_\alpha)_{A^\Gamma}}\Ms\ |\ \Gamma_{(\As_\beta)_{B^\Gamma}}\xrightarrow{\tilde{T}} {}_{(\As_\alpha)_{A^\Gamma}}\Ns_{(\As_\beta)_{B^\Gamma}}$ of multimodules such that $T=\tilde{T}\circ T^\Gamma_\Ms$.
\end{definition}

\begin{remark} \label{rem: contractions}
As usual, $\Gamma$-contractions of $\Zs$-central multimodules are unique up to a unique isomorphisms compatible with the defining factorization property; their existence is provided by the following construction. 

\medskip 

Given the $\Zs$-central multimodule ${}_{(\As_\alpha)_{A}}\Ms_{(\As_\beta)_{B}}$ and the injective symmetric relation $\Gamma\subset (A\uplus B)\times(A\uplus B)$ with $\As_\xi=\As_\zeta$ whenever $(\xi,\zeta)\in\Gamma$, defining $A^\Gamma:=A-\dom(\Gamma)$ and $B^\Gamma:=B-\im(\Gamma)$, consider the $\Zs$-central 
$(\As_\alpha)_{\alpha\in A^\Gamma}$-$(\As_\beta)_{\beta\in B^\Gamma}$ sub-multimodule ${}_{(\As_\alpha)_{A^\Gamma}}[\Ms,\Gamma]_{(\As_\beta)_{B^\Gamma}}$ of ${}_{(\As_\alpha)_{A^\Gamma}}\Ms_{(\As_\beta)_{B^\Gamma}}$ generated by the elements of the form:
\begin{align*}
a\cdot_\xi x&-a\cdot_\zeta x, \quad \forall x\in\Ms, \ \forall a\in \As_\xi=\As_\zeta, \quad 
(\xi,\zeta)\in\Gamma\cap (A\times A), 
\\
x\cdot_\xi a&-x\cdot_\zeta a, \quad \forall x\in\Ms, \ \forall a\in \As_\xi=\As_\zeta, \quad 
(\xi,\zeta)\in\Gamma\cap (B\times B),  
\\
a\cdot_\xi x&-x\cdot_\zeta a, \quad \forall x\in\Ms, \ \forall a\in \As_\xi=\As_\zeta, \quad 
(\xi,\zeta)\in\Gamma\cap (A\times B). 
\end{align*}
The quotient map $T^\Gamma_\Ms:\Ms\to \frac{\Ms}{[\Ms,\Gamma]}=:\Ms\ | \ \Gamma$ onto the quotient $(\As_\alpha)_{A^\Gamma}$-$(\As_\beta)_{B^\Gamma}$ multimodule, satisfies the universal factorization property, since every $\Gamma$-tracial homomorphism $\Ms\xrightarrow{T}\Ns$ of  $(\As_\alpha)_{A^\Gamma}$-$(\As_\beta)_{B^\Gamma}$ multimodules entails $[\Ms,\Gamma]\subset \ke(T)$ and hence canonically factorizes via $T^\Gamma_\Ms$. 

\medskip 

Notice that it is possible to have multimodules that only possess trivial $\Gamma$-traces (for a certain family $\Gamma$) and hence they have trivial universal $\Gamma$-contractions. 
\xqed{\lrcorner}
\end{remark}

We briefly examine how involutions in multimodules descend to their contractions. 
\begin{definition}
An involution $(\dagger^\sigma_\Ms,\ddagger^\sigma_\Ms,\star_\Ms)_{\sigma\in f_\Ms}$ on a $\Zs$-central multimodule ${}_{(\As_\alpha)_{A}}\Ms_{(\As_\beta)_{B}}$ over $\Rs_\Zs$-algebras is a \emph{$\Gamma$-compatible involution}, where 
$\Gamma\subset (A\uplus B)\times(A\uplus B)$ is an injective symmetric relation 
on $A\uplus B$, if: 
\begin{equation} \label{eq: Gamma} 
(\xi,\zeta)\in\Gamma \imp (f(\xi),f(\zeta))\in\Gamma, \quad \forall \xi,\zeta\in A\uplus B. 
\end{equation}
\end{definition}

\begin{proposition}\label{prop: inv-tr}
Suppose that $({}_{(\As_\alpha)_{A}}\Ms_{(\As_\beta)_{B}},\star_\Ms)$ in an involutive $\Zs$-central multimodule over $\Rs_\Zs$-algebras. If the involution is $\Gamma$-compatible with a $\Gamma$-contraction, there exists a unique contracted involution onto $\Ms\ |\ \Gamma$ and the contraction map $\Ms\xrightarrow{T_\Gamma} \Ms\ | \ \Gamma$ is involutive.  
\end{proposition}
\begin{proof}
Condition~\eqref{eq: Gamma} implies that the involution $\star_\Ms$ leaves invariant the submultimodule ${}_{{(\As_\alpha)}_{A^\Gamma}}[\Ms,\Gamma]_{{(\Bs_\beta)}_{B^\Gamma}}$ and hence, defining $(x+[\Ms,\Gamma])^{\star}:=x^{\star_\Ms}+[\Ms,\Gamma]$, for all $x\in\Ms$, the involution will pass to the quotient multimodule (with the same covariance properties in $(\alpha,\beta)\in A^\Gamma\times B^\Gamma$) and $T_\Gamma(x^{\star_\Ms})=T_\Gamma(x)^\star$, for all $x\in\Ms$. 
\end{proof}

\begin{remark}\label{rem: inv-tr}
Looking at the universal contructions of tensor products of multimodules in remark~\ref{rem: otimes-tr} and of contrations in remark~\ref{rem: contractions}, we obtain the following familiar result: 
\begin{equation*}
\left({}_{(\As_\alpha)_{\alpha\in A}}\Ms_{(\Bs_\beta)_{\beta\in B}} \right)
\otimes_\Gamma 
\left({}_{(\Cs_\gamma)_{\gamma\in C}}\Ms_{(\Ds_\delta)_{\delta\in D}} \right)
\simeq 
T_\Gamma \left({}_{(\As_\alpha)_{\alpha\in A}}\Ms_{(\Bs_\beta)_{\beta\in B}} \otimes_\Zs 
{}_{(\Cs_\gamma)_{\gamma\in C}}\Ms_{(\Ds_\delta)_{\delta\in D}} \right), 
\end{equation*}
tensor products over $\Gamma\subset (A\uplus C)\uplus(B\uplus D)$ are naturally isomorphic to $\Gamma$-contracted tensor products over $\Zs$.
\xqed{\lrcorner}
\end{remark}

\bigskip 

We pass now to consider the generalization of inner product couplings for multimodules. 

\medskip 

\begin{definition} \label{def: ip-multi}
Suppose that the unital associative $\Rs_\Zs$-algebras $\As_\alpha$ and $\Bs_\beta$, for $(\alpha,\beta)\in A\times B$, are all contravariantly involutive.  
Given a multimodule ${}_{(\As_\alpha)_A}\Ms_{(\Bs_\beta)_B}$ and sub-indexes $I\times J\subset A\times B$, an \emph{$(I,J)$-right-inner product on $\Ms$} consists of a bi-additive map $\ip{\cdot}{\cdot}_{I-J}: \Ms\times\Ms\to \left(\bigotimes^\Zs_{i\in I}\As_i\right)\otimes_\Zs\left(\bigotimes^\Zs_{j\in J}\Bs_j\right)$ such that: 
\begin{align*}
\ip{x}{a\cdot_\alpha y\cdot_\beta b}_{I-J} &= a\cdot_\alpha\ip{x}{y}_{I-J}\cdot_\beta b,  
&&
\forall x,y\in\Ms, \ (a,b)\in\As_\alpha\times\Bs_\beta, \ (\alpha,\beta)\in I\times J, 
\\
\ip{a\cdot_\alpha x\cdot_\beta b}{y}_{I-J} &=b^*\bullet_\beta\ip{x}{y}_{I-J}\bullet_\alpha a^*, 
&& 
\forall x,y\in\Ms, \ (a,b)\in\As_\alpha\times\Bs_\beta, \ (\alpha,\beta)\in I\times J, 
\\
\ip{a\cdot_\alpha x\cdot_\beta b}{y}_{I-J} &=\ip{x}{a^*\cdot_\alpha y\cdot_\beta b^*}_{I-J},
&& 
\forall x,y\in\Ms, \ (a,b)\in\As_\alpha\times\Bs_\beta, \ (\alpha,\beta)\in (A-I)\times (B-J).  
\end{align*}
A \emph{$(I,J)$-left-inner product on $\Ms$} is a bi-additive map ${}_{I-J}\ip{\cdot}{\cdot}: \Ms\times\Ms\to \left(\bigotimes^\Zs_{i\in I}\As_i\right)\otimes_\Zs\left(\bigotimes^\Zs_{j\in J}\Bs_j\right)$ such that: 
\begin{align*}
{}_{I-J}\ip{a\cdot_\alpha x\cdot_\beta b}{y} &= a\cdot_\alpha{}_{I-J}\ip{x}{y}\cdot_\beta b,  
&&
\forall x,y\in\Ms, \ (a,b)\in\As_\alpha\times\Bs_\beta, \ (\alpha,\beta)\in I\times J, 
\\
{}_{I-J}\ip{x}{a\cdot_\alpha y\cdot_\beta b} &=b^*\bullet_\beta{}_{I-J}\ip{x}{y}\bullet_\alpha a^*, 
&& 
\forall x,y\in\Ms, \ (a,b)\in\As_\alpha\times\Bs_\beta, \ (\alpha,\beta)\in I\times J, 
\\
{}_{I-J}\ip{a\cdot_\alpha x\cdot_\beta b}{y} &={}_{I-J}\ip{x}{a^*\cdot_\alpha y\cdot_\beta b^*},
&& 
\forall x,y\in\Ms, \ (a,b)\in\As_\alpha\times\Bs_\beta, \ (\alpha,\beta)\in (A-I)\times (B-J).  
\end{align*}
\end{definition}

\begin{proposition}
For every right-$(I,J)$-inner product $(x,y)\mapsto\ip{x}{y}_{I-J}$ we have its:
\begin{equation*} 
\text{\emph{transpose}} \quad (x,y)\mapsto\ip{y}{x}_{I-J}, 
\quad
\text{\emph{$*$-conjugate}} \quad (x,y)\mapsto\ip{x}{y}^*_{I-J}, 
\quad
\text{\emph{$*$-adjoint}} \quad (x,y)\mapsto\ip{y}{x}^*_{I-J}.  
\end{equation*}
The transpose and conjugate are left-$(I,J)$-inner products; the adjoint is a right-$(I,J)$-inner product on $\Ms$.
\end{proposition}

\begin{remark}
Without entering into a detailed discussion of ``positivity'' for inner products, we simply mention that stating this condition requires an additional compatible ``order structure'' on the involved rings and algebras. Whenever the commutative unital associative involutive ring $\Zs$ is equipped with a positive cone $\Zs_+$ (that by definition is a pointed subset $0_\Zs\in \Zs_+\subset \Zs$, stable under addition $\Zs_+ + \Zs_+\subset\Zs_+$, stable under multiplication $\Zs_+\cdot\Zs_+\subset\Zs_+$, sharp $\Zs_+\cap(-\Zs_+)=\{0_\Zs\}$ and involutive $\Zs_+^{\star_\Zs}\subset\Zs_+$), any $\Zs$-central unital associative algebra $\Rs$ (and hence any $\Zs$-central $\Rs_\Zs$-algebra) canonically inherits a positive cone $\Rs_+:=\Zs_+\cdot 1_\Rs\subset \Rs$.  
In this case, positivity of a right-$(I,J)$-inner product can be imposed requiring $\ip{x}{x}_{I-J}\in\left[ \left(\bigotimes^\Zs_{i\in I}\As_i\right)\otimes_\Zs\left(\bigotimes^\Zs_{j\in J}\Bs_j\right)\right]_+$, for all $x\in\Ms$. Similar condition can be imposed for left-$(I,J)$-inner products. 
\xqed{\lrcorner}
\end{remark}

\begin{theorem}\label{th: Riesz}
An $(I,J)$-inner product (right or left) induces canonical \emph{Riesz maps} 
\begin{gather*}
\Ms\xrightarrow{{}^I\overrightarrow{\Lambda}^J} {}^{*_I}\Ms^{*_J},
\quad \quad {}^I\overrightarrow{\Lambda}_x^J:x\mapsto {}^I\overrightarrow{\Lambda}_x^J\quad \quad  {}^I\overrightarrow{\Lambda}_x^J:y\mapsto \ip{x}{y}_{I-J}, 
\\ 
\Ms\xrightarrow{{}^I\overleftarrow{\Lambda}^J} {}^{*_I}\tilde{\Ms}^{*_J},
\quad \quad  {}^I\overleftarrow{\Lambda}_y^J:y\mapsto {}^I\overleftarrow{\Lambda}_y^J\quad \quad  {}^I\overleftarrow{\Lambda}_y^J:x\mapsto \ip{x}{y}_{I-J}, 
\end{gather*}
where ${}^I\overrightarrow{\Lambda}^J$ is a contravariant morphism of multimodules into ${}^{*_I}\Ms^{*_J}$, the $(I,J)$-dual of $\Ms$, and respectively ${}^I\overleftarrow{\Lambda}^J$ is a covariant morphism of multimodules into the ${}^I\overrightarrow{\Lambda}^J$-twisted of ${}^{*_I}{\Ms}^{*_J}$, here denoted by ${}^{*_I}\tilde{\Ms}^{*_J}$.  
\end{theorem}
\begin{proof}
From definition~\ref{def: ip-multi} we see that ${}^I\overrightarrow{\Lambda}^J_x\in {}^{*_I}\Ms^{*_J}$, for all $x\in\Ms$.  The map ${}^I\overrightarrow{\Lambda}^J:\Ms\to {}^{*_I}\Ms^{*_J}$ is a contravariant morphism of multimodules: 
\begin{align*}
{}^I\overrightarrow{\Lambda}^J_{a\cdot_\alpha x\cdot_\beta b}(y) 
&=\ip{a\cdot_\alpha x\cdot_\beta b}{y}_{I-J}
=b^*\bullet_\beta \ip{x}{y}_{I-J}\bullet_\alpha a^*
=b^*\bullet_\beta {}^I\overrightarrow{\Lambda}^J_x (y)\bullet_\alpha a^*
\\
&=(b^*\bullet_\beta {}^I\overrightarrow{\Lambda}^J_x \bullet_\alpha a^*)(y), \quad 
\forall (\alpha,\beta)\in I\times J, \ x,y\in\Ms, \ (a_\alpha,b_\beta)\in\As_\alpha\times \Bs_\beta; 
\\
{}^I\overrightarrow{\Lambda}^J_{a\cdot_\alpha x\cdot_\beta b}(y) 
&=\ip{a\cdot_\alpha x\cdot_\beta b}{y}_{I-J}
=\ip{x}{a^*\cdot_\alpha y\cdot_\beta b^*}_{I-J}
={}^I\overrightarrow{\Lambda}^J_x(a^*\cdot_\alpha y\cdot_\beta b^*)
\\
&=(b^*\odot_\beta {}^I\overrightarrow{\Lambda}^J_x \odot_\alpha a^*)(y), \quad
\forall (\alpha,\beta)\in (A-I)\times (B-J), \ x,y\in\Ms, \ (a_\alpha,b_\beta)\in\As_\alpha\times \Bs_\beta.
\end{align*}
The proof for the case of ${}^I\overleftarrow{\Lambda}^J$ can be obtained passing to the transpose inner product. 
\end{proof}

\begin{definition}
An inner product $\ip{\cdot}{\cdot}_{I-J}$ is \emph{$*$-Hermitian} if it coincides with its $*$-adjoint; \emph{non-degenerate} if both the Riesz maps ${}^I\overrightarrow{\Lambda}^J$ and ${}^I\overleftarrow{\Lambda}^J$ are injective; 
\emph{algebraically full} if $\ip{\Ms}{\Ms}_{I-J}=\left(\bigotimes^\Zs_{i\in I}\As_i\right)\otimes_\Zs\left(\bigotimes^\Zs_{j\in J}\Bs_j\right)$; \emph{saturated} if both 
${}^I\overrightarrow{\Lambda}^J$ and ${}^I\overleftarrow{\Lambda}^J$ are surjective. 
\end{definition}

\begin{remark}
Notice that the Riesz map ${}^I\overrightarrow{\Lambda}^J$ is contravariant and hence, under non-degeneracy and saturation, an $(I,J)$-inner product always induces an anti-isomorphism between $\Ms$ and its $(I,J)$ dual ${}^{*_I}\Ms^{*_J}$. 

Under fullness and saturation, ${}^{*_I}\Ms^{*_J}\xrightarrow{({}^I\overrightarrow{\Lambda}^J)^{-1}}\Ms$ is a  $({}^I\overrightarrow{\Lambda}^J)^{-1}$-twisted of $\Ms$ as defined in~\ref{def: twisted} and~\ref{prop: cj-duals}. 

\medskip 

The contravariant nature of Riesz maps requires contravariant involutions in the definition of inner products;
alternative possibilities can be explored with ``inner couplings'' on $\Ms$ with more general signatures.
\xqed{\lrcorner}
\end{remark}

\begin{remark}\label{rem: hybrid}
Thinking of multimodules in the 1-category ${}_{(\As_\alpha)_A}\Mf_{(\Bs_\beta)_B}$ as ``1-arrows'' $(\As_\alpha)_A\xleftarrow{\Ms}(\Bs_\beta)_B$, in a 2-category of morphisms $\xymatrix{(\As_\alpha)_A & \ltwocell^{\Ms}_{\Ns}{\Phi} (\Bs_\beta)_B}$, we see that, for all $I\times J\in A\times B$, $(I,J)$-duals provide an $\{0,1\}$-contravariant involution $\xymatrix{(\As_\alpha)_A  \rrtwocell^{{}^{*_I}\Ns^{*_J}}_{{}^{*_I}\Ms^{*_J}}{\quad {}^{*_I}\Phi^{*_J}} & & (\Bs_\beta)_B}$, over objects and 1-arrows, in the sense described in~\cite[section~4]{BCLS20}. 
Riesz maps can be considered as ``natural transformations'' examples of \textit{hybrid 2-arrows} 
$\xymatrix{(\As_\alpha)_A  \rrtwocell^{\Ms}_{{}^{*_I}\Ms^{*_J}}{'\ \quad {\Downarrow}\  {}^{*_I}\overrightarrow{\Lambda}^{*_J}} & & (\Bs_\beta)_B}$
following the definition of \textit{hybrid 2-category} described in~\cite{BePu14}. We will pursue such developments elsewhere.
\xqed{\lrcorner}
\end{remark}

\begin{remark}
In definition~\ref{def: ip-multi} in order to keep the closest possible resemblance to the usual axioms for inner products in Hilbert spaces and Hilbert-C*-modules, we have imposed covariance, for certain actions, only one of the two variables and contravariance on the other. It is perfectly possible to consider more general cases, where covariance and contravariance are simultaneously present in both variables (on disjoint sets of indexes): 
let $I\times J\subset A\times B$ with $I:=I_l\cup I_r$, $J=J_l\cup J_r$ and 
$I_l\cap I_r=\varnothing=J_l\cap J_r$,  
a \emph{$(I_l,J_l)$-left $(I_r,J_r)$-right inner product} on ${}_{(\As_\alpha)_A}\Ms_{(\Bs_\beta)_B}$ is a bi-additive map $(x,y)\mapsto {}_{I_l-J_l}\ip{x}{y}_{I_r-J_r}$, for $x,y\in \Ms$, such that:  $\forall x,y\in\Ms,$
\begin{align*}
{}_{I_l-J_l}\ip{x}{a\cdot_\alpha y\cdot_\beta b}_{I_r-J_r} &= a\cdot_\alpha\ip{x}{y}_{I_r-J_r}\cdot_\beta b,  
&&
\forall (a,b)\in\As_\alpha\times\Bs_\beta, \ (\alpha,\beta)\in I_r\times J_r, 
\\
{}_{I_l-J_l} \ip{a\cdot_\alpha x\cdot_\beta b}{y}_{I_r-J_r}&= a\cdot_\alpha\ip{x}{y}_{I_r-J_r}\cdot_\beta b,  
&&
\forall (a,b)\in\As_\alpha\times\Bs_\beta, \ (\alpha,\beta)\in I_l\times J_l, 
\\
{}_{I_l-J_l}\ip{a\cdot_\alpha x\cdot_\beta b}{y}_{I_r-J_r} &=b^*\bullet_\beta{}_{I_l-J_l}\ip{x}{y}_{I_r-J_r}\bullet_\alpha a^*, 
&& 
\forall (a,b)\in\As_\alpha\times\Bs_\beta, \ (\alpha,\beta)\in I_r\times J_r, 
\\ 
{}_{I_l-J_l}\ip{x}{a\cdot_\alpha y\cdot_\beta b}_{I_r-J_r} &=b^*\bullet_\beta{}_{I_l-J_l}\ip{x}{y}_{I_r-J_r}\bullet_\alpha a^*, 
&& 
\forall (a,b)\in\As_\alpha\times\Bs_\beta, \ (\alpha,\beta)\in I_l\times J_l, 
\\
{}_{I_l-J_l}\ip{a\cdot_\alpha x\cdot_\beta b}{y}_{I_r-J_r} &={}_{I_l-J_l}\ip{x}{a^*\cdot_\alpha y\cdot_\beta b^*}_{I_r-J_r},
&& 
\forall (a,b)\in\As_\alpha\times\Bs_\beta, \ (\alpha,\beta)\in (A-I)\times (B-J).  
\end{align*}
Riesz maps ${}^{I_l-J_l}\overrightarrow{\Lambda}^{I_r-J_r}$ and ${}^{I_l-J_l}\overleftarrow{\Lambda}^{I_r-J_r}$ can be similarly defined and a perfect parallel of theorem~\ref{th: Riesz} holds. 
\xqed{\lrcorner}
\end{remark}

\section{Outlook} \label{sec: outlook}

Although we are not going here into specific details, that will be subject of a forthcoming work, we preview some of the categorical features making multimodules a quite intriguing playground. 

\medskip  

The family of multimodules, with their several tensor products, constitutes a paradigmatic example of ``algebraic structure'' consisting of ``many inputs / many outputs nodes'' that can be ``linked'' in many different ways: each multimodule ${}_{(\As_\alpha)_A}\Ms_{(\Bs_\beta)_B}$ should be interpreted as an arrow with sources $(\Bs_\beta)_B$ and targets $(\As_\alpha)_A$;\footnote{
We are using here, for the tensor products, the same ``reversed order'' notation of the functional compositions.
} 
every tensor product over a subfamily provides a possible ``concatenation'' of arrows and such compositions will be subject to associativity and unitality axioms typical of category theory. 

\medskip 

At the 1-categorical level (when only multimodules as 1-arrows and tensor products as compositions are considered) the structure seems to be describable as a \textit{colored properad}~\cite{HRY15}, a horizontal categorification (i.e.~a many objects version) of the notion of \textit{properad} introduced by~\cite{Va07}. 

\medskip 

Dualities of multimodules seem to provide the easiest examples of involutions for arrows in a colored properad and can be taken as a paradigmantic template in order to axiomatize a notion of ``involutive colored properad''. 
Contractions can be used to introduce ``sinks'' and ``sources'', hence more general types of ``partial involutions''.  

\medskip 

As already mentioned in remark~\ref{rem: hybrid}, we plan to further study Riesz dualities as examples of \textit{hybrid natural transformations}, between functors with different covariance, in the context of hybrid \hbox{2-categories} introduced in~\cite{BePu14}. 

\medskip 

Covariant morphisms of multimodules should be interpreted (exactly as in the usual case of categories of bimodules) as cubical 2-arrows. In this way, one obtains for multimodules a colored properad analog of the usual double category of covariant morphisms of bimodules. 

\medskip 

It is also possible to iterate the construction of multimodules over multimodules creating a vertical categorification ladder that can be used to define ``involutive higher colored properads'' (possibly requiring the non-commutative exchange property introduced in~\cite{BCLS20}).  

\medskip 

The purely algebraic theory of $\Zs$-central multimodules over $\Rs_\Zs$-algebras here presented can be subject to a functional analytic treatment as soon as topologies/uniformities are introduced and the actions are required to be continuous in the suitable sense. We will explore in the future the more restrictive axioms for a (higher) C*-algebraic version of this material and obtain (infinite-dimensional) functional analytic generalizations of the (essentially finite-dimensional) \textit{reflexivity} ($\theta^\Ms$ covariant isomorphism) and \textit{self-duality} (${}^J\overrightarrow{\Lambda}^I_\Ms$ contravariant isomorphism) conditions on multimodules. 

\medskip 

Although the basic definition of first-order differential operator between $\Zs$-central multimodules over non-commutative $\Rs_\Zs$-algebras is included in appendix~\ref{sec: 1st-ord-multi}, much more needs to be done regarding the full differential theoretic theory of multimodules (and also bimodules!), starting with a theory of connections on multimodules and possibly proceeding in the direction of \textit{properadic non-commutative geometry} as a natural extension of our current efforts in categorical non-commutative geometry. An exploration of the interplay between duality (for bimodules) and first-order differential operators associated to covariant differential calculi on a non-commutative $\Zs$-central algebra is carried on in our forthcoming work (mentioned in footnote~\ref{foo: ncdc}). 

\bigskip 


\emph{Notes and Acknowledgments:} 
%
%
%
%
%
P.Bertozzini thanks Starbucks Coffee ($1^{\text{st}}$ floor of Emporium Tower, Emquartier Sky Garden, Jasmine City) where he spent most of the time dedicated to this research project; he thanks Fiorentino Conte of ``The Melting Clock'' for the great hospitality during many crucial on-line dinner-time meetings Bangkok-Rome. 



\appendix

\section{Functorial Pairings and Semi-adjunctions} \label{sec: functorial-pairing}

In order to properly discuss the categorical features of dual pairing of multimodules, we need to deal with a variant of the well-known notion of adjunction between functors, originally introduced in~\cite{Med74}, and generalized (in a much wider context) in~\cite{MeWi13,Wis13}.  

\begin{definition}
A \emph{full functorial pairing} $\Fg \stackrel[\rho]{\lambda}{\leftrightarrows}|\ \Gg$ between the covariant functors $\xymatrix{\Af\rtwocell^{\Fg}_{\Gg}{'} & \Bf }$ is a pair of natural transformations between the (left-contravariant, right-covariant) $\Hom$-bifunctors~\cite{Wis13}: 
\begin{equation}\label{eq: semi-adj}
\Hom_\Bf(\Fg(A),B)\stackrel[\rho_{AB}]{\lambda_{AB}}{\leftrightarrows}\Hom_\Af(A,\Gg(B)), 
\quad \forall (A,B)\in\Ob_\Af\times\Ob_\Bf.
\end{equation}
The full pairing is \emph{regular} if both $\rho,\lambda$ are regular maps: 
\begin{equation*}
\rho_{AB}\circ\lambda_{AB}\circ\rho_{AB}=\rho_{AB},  
\quad \lambda_{AB}\circ\rho_{AB}\circ\lambda_{AB}=\lambda_{AB}, 
\quad \forall (A,B)\in\Ob_\Af\times\Ob_\Bf.
\end{equation*}
A covariant full functorial pairing $\Fg \stackrel[\rho]{\lambda}{\leftrightarrows}|\ \Gg$ will be called~\cite{Med74} a: 
\begin{align*}
&\text{\emph{right semi-adjunction $\Fg \stackrel[\rho]{\lambda}{\leftrightarrows}|\ \underline{\Gg}$} if:} 
&\rho_{AB}\circ\lambda_{AB}=\id_{\Hom_\Af(A,\Gg(B))}, \quad  \forall (A,B)\in\Ob_\Af\times\Ob_\Bf, 
\\
&\text{\emph{left semi-adjunction $\underline{\Fg} \stackrel[\rho]{\lambda}{\leftrightarrows}|\ \Gg$} if:} 
&\lambda_{AB}\circ\rho_{AB}=\id_{\Hom_\Bf(\Fg(A),B)}, \quad  \forall (A,B)\in\Ob_\Af\times\Ob_\Bf.
\end{align*}
\end{definition}

The following remark follows immediately from~\cite[section~2]{Wis13}. 
\begin{remark}
Note that, similarly to adjunction (as indicated by our notation $\Fg \stackrel[\rho]{\lambda}{\leftrightarrows}|\ \Gg$),  full functorial pairing is an ``asymmetric'' notion, with the functor $\Gg$ right functorially paired to $\Fg$ (equivalently $\Fg$ left paired to $\Gg$). 

\medskip 

The existence of the natural transformation $\rho$ (respectively $\lambda$) in formula~\eqref{eq: semi-adj} is equivalent to the existence of a unit $\id_{\Af}\xrightarrow{\eta}\Gg\circ\Fg$ (respectively a co-unit $\Fg\circ\Gg\xrightarrow{\epsilon}\id_{\Bf}$) natural transformation: 
given the full functorial pairing $\Fg \stackrel[\rho]{\lambda}{\leftrightarrows}|\ \Gg$, one defines, for all $(A,B)\in\Ob_\Af\times\Ob_\Bf$,  
\begin{align*}
&\eta_A:=\rho_{A\Fg(A)}(\iota_{\Fg(A)}), 
&
\epsilon_B:=\lambda_{\Gg(B)B}(\iota_{\Gg(B)}); 
\end{align*}
in the reverse direction, given unit  $\id_{\Af}\xrightarrow{\eta}\Gg\circ\Fg$ and co-unit $\Fg\circ\Gg\xrightarrow{\epsilon}\id_{\Bf}$ natural transformations, one defines
\begin{align*}  
&\rho_{AB}(x):=\Gg(x)\circ\eta_A\in\Hom_\As(A;\Gg(B)), \quad \forall (A,B)\in\Ob_\Af\times\Ob_\Bf, \ \forall x\in\Hom_\Bf(\Fg(A);B), 
\\
&\lambda_{AB}(y):=\epsilon_B\circ \Fg(y)\in\Hom_\Bf(\Fg(A);B), \quad \forall (A,B)\in\Ob_\Af\times\Ob_\Bf, 
 \ \forall y\in\Hom_\Af(A;\Gg(B)). 
\end{align*}

\medskip 

The semi-adjunction conditions, can be equivalently written via composition of units and co-units in the \hbox{2-category} of natural transformations as:\footnote{
These are actually the original equations used by~\cite{Med74} to define right and left semi-adjunctions.
}
\begin{align*}
&
\Fg\ \stackrel[\rho]{\lambda}{\leftrightarrows}|\ \Gg  
\quad A\xrightarrow{\eta_A:=\rho_{A\Fg(A)}(\iota_{\Fg(A)})}\Gg(\Fg(A)), 
&&
B\xleftarrow{\epsilon_B:=\lambda_{\Gg(B)B}(\iota_{\Gg(B)})}\Fg(\Gg(B)),   
\\
& \underline{\Fg}\ \stackrel[\rho]{\lambda}{\leftrightarrows} | \ \Gg \quad \text{$\Fg$ is left semi-adjoint of $\Gg$:} 
&&\epsilon_{\Fg(A)}\circ_\Bf \Fg(\eta_A)=\iota_{\Fg(A)} \iff \lambda\circ\rho=\id_{\Hom_\Bf(\Fg(A);B)}; 
\\
&
\Fg\ \stackrel[\rho]{\lambda}{\leftrightarrows} | \ \underline{\Gg} \quad
\text{$\Gg$ is right semi-adjoint of $\Fg$:}
&&
\Gg(\epsilon_B)\circ_\Af \eta_{\Gg(B)}=\iota_{\Gg(B)} \iff \rho\circ \lambda=\id_{\Hom_\Af(A;\Gg(B))}. 
\end{align*}
In view of the established equivalence description of semi-adjunctions via $(\lambda,\rho)$ or via $(\eta,\epsilon)$ we will liberally utilize the alternative notations: 
\begin{equation*}
\Fg \stackrel[\rho]{\lambda}{\leftrightarrows}|\ \Gg \ \iff:\ \Fg \stackrel[\eta]{\epsilon}{\leftrightarrows}|\ \Gg
\end{equation*}

\medskip 

A semi-adjunction is necessarily a regular full functorial pairing; whenever $\rho$ and $\lambda$ are inverse of each other, the necessarily regular full functorial pairing reproduces an adjunction $\Fg\dashv \Gg$ with unit $\eta$ and co-unit $\epsilon$.
\xqed{\lrcorner}
\end{remark}

We will need to utilize semi-adjunctions in the case of contravariant functors. 
\begin{remark}\label{rem: s-aj}
In the case of contravariant functors $\xymatrix{\Af\rtwocell^{\Fg}_{\Gg}{'} & \Bf}$, the usual right-right and left-left adjunctions, are corresponding to a \emph{right contravariant full functorial pairing}: 
\begin{equation*} 
\Fg\ | \stackrel[\rho]{\lambda}{\leftrightarrows}|\ \Gg, 
\quad \quad 
\Hom_\Bf(B,\Fg(A))\stackrel[\rho_{AB}]{\lambda_{AB}}{\leftrightarrows}\Hom_\Af(A,\Gg(B)), 
\quad\forall (A,B)\in\Ob_\Af\times\Ob_\Bf, 
\end{equation*}
and respectively to a \emph{left contravariant functorial pairing}: 
\begin{equation*} 
\Fg\  \stackrel[\rho]{\lambda}{\leftrightarrows}\ \Gg, 
\quad \quad 
\Hom_\Bf(\Fg(A),B)\stackrel[\rho_{AB}]{\lambda_{AB}}{\leftrightarrows}\Hom_\Af(\Gg(B),A),
\quad\forall (A,B)\in\Ob_\Af\times\Ob_\Bf. 
\end{equation*}
The definitions of contravariant regularity and contravariant semi-adjunction remain the same. 

\medskip 

For all possible cases of contravariant semi-adjunction, the equivalent statements in terms of the associated unit and co-unit natural transformations can be summarized as follows, for all $(A,B)\in\Ob_\Af\times\Ob_\Bf$: 
\begin{align}
\notag
& 
\Fg\ | \stackrel[\eta]{\epsilon}{\leftrightarrows}|\ \Gg  
\quad \quad 
A\xrightarrow{\eta_A:=\rho_{A\Fg(A)}(\iota_{\Fg(A)})}\Gg(\Fg(A)), 
&&
B\xrightarrow{\epsilon_B:=\lambda_{\Gg(B)B}(\iota_{\Gg(B)})}\Fg(\Gg(B)), 
\\ \label{eq: r-s-aj}
&\underline{\Fg}\ | \stackrel[\eta]{\epsilon}{\leftrightarrows}|\ \Gg 
\quad 
\text{$\Fg$ is right semi-adjoint of $\Gg$:} 
&&
\Fg(\eta_A)\circ_\Bf\epsilon_{\Fg(A)} =\iota_{\Fg(A)} \iff \lambda_{AB}\circ\rho_{AB}=\id_{\Hom_\Bf(\Fg(A);B)},  
\\ \notag
&\Fg \ | \stackrel[\eta]{\epsilon}{\leftrightarrows}|\ \underline{\Gg}
\quad \text{$\Gg$ is right semi-adjoint of $\Fg$:}
&&
\Gg(\epsilon_B)\circ_\Af\eta_{\Gg(B)}=\iota_{\Gg(B)} \iff \rho_{AB}\circ \lambda_{AB}=\id_{\Hom_\Af(A;\Gg(B))};  
\\ \notag 
&\Fg\  \stackrel[\eta]{\epsilon}{\leftrightarrows} \ \Gg  
\quad \quad 
A\xleftarrow{\eta_A:=\rho_{A\Fg(A)}(\iota_{\Fg(A)})}\Gg(\Fg(A)), 
&&
B\xleftarrow{\epsilon_B:=\lambda_{\Gg(B)B}(\iota_{\Gg(B)})}\Fg(\Gg(B)), 
\\ \notag
&\Fg \ \stackrel[\eta]{\epsilon}{\leftrightarrows}\ \underline{\Gg}
\quad \text{$\Gg$ is left semi-adjoint of $\Gg$:}
&&
\epsilon_{\Fg(A)}\circ_\Bf \Fg(\eta_A)=\iota_{\Fg(A)} \iff \lambda_{AB}\circ\rho_{AB}=\id_{\Hom_\Bf(\Fg(A);B)},
\\ \notag 
&\underline{\Fg} \ \stackrel[\eta]{\epsilon}{\leftrightarrows}\ \Gg
\quad \text{$\Fg$ is left semi-adjoint of $\Fg$:}
&& 
\eta_{\Gg(B)}\circ_\Af\Gg(\epsilon_B)=\iota_{\Gg(B)} \iff \rho_{AB}\circ \lambda_{AB}=\id_{\Hom_\Af(A;\Gg(B))}.  
\end{align}
Notice, due to the contravariance, the ``change of direction'' and respectively the order of composition of the natural transformations involved (so that for contravariant right functorial pairings we have in practice two unit, and for contravariant left functorial pairings, we actually have two co-unit natural transformations).

\medskip 

Notice also that the ``doubling'' of the semi-adjointness conditions is just an apparent artifact of notation since: 
$\underline{\Fg}\ | \stackrel[\rho]{\lambda}{\leftrightarrows}|\ \Gg \iff 
\underline{\Gg}\ | \stackrel[\lambda]{\rho}{\leftrightarrows}|\ \Fg $ and similarly 
$\underline{\Fg}\ \stackrel[\rho]{\lambda}{\leftrightarrows} \ \Gg \iff 
\underline{\Gg}\  \stackrel[\lambda]{\rho}{\leftrightarrows} \ \Fg$. 

\medskip 

Whenever $\rho$ and $\lambda$ are inverse of each other, the necessarily regular full contravariant right functorial pairing reproduces a contravariant right adjunction $\Fg\vdash \dashv \Gg$ with two units $\eta$ and $\epsilon$ that, upon restriction to the full subcategories of reflexive objects provides a duality. 
\xqed{\lrcorner}
\end{remark}

\medskip 

Let us exemplify the required semi-adjunction in the case of $\Zs$-central bimodules over $\Zs$-central $\Rs$-algebras.
\begin{proposition}\label{prop: sadj-z}
In the category $\Mf_{\Zs}$ of morphisms of $\Zs$-central $\Zs$-bimodules, the transposition pairing duality $\Omega\mapsto \Omega^*$, for $\Omega\in\Ob(\Mf_{\Zs})$, is a contravariant endofunctor. The evaluation natural transformation $\ev$ induces a right contravariant semi-adjoint endo-funtorial pairing $\xymatrix{\Mf_\Zs \rtwocell^{{}_*}_{{}^*}{'} & \Mf_\Zs}$, hence the transposition duality endofunctor is right semi-adjoint to itself:
$\underline{*} \ | \stackrel[\ev]{\ev}{\leftrightarrows}\ |\ * \iff  *\ | \stackrel[\ev]{\ev}{\leftrightarrows}\ |\ \underline{*}$.
Upon restriction to the subcategory of reflexive objects (those objects $\Omega$ for which $\ev^\Omega$ is an isomorphism) the functorial pairing above is a duality. 
\end{proposition}
\begin{proof}
The transposition duality functor associates to every morphism $\Omega_1\xrightarrow{\Phi}\Omega_2$ of $\Zs$-central $\Zs$-bimodules the transposed map $\Omega_1^*\xleftarrow{\Phi^*}\Omega_2^*$ defined, for all $\psi\in\Omega_2^*$, by $\Phi^*(\psi):=\psi\circ\Phi\in\Omega_1^*$; the map $\Phi^*$ is $\Zs$-linear and the transposition $\Phi\mapsto \Phi^*$ is a contravariant endofunctor: $(\Phi\circ\Psi)^*=\Psi^*\circ\Phi^*$ and $(\id_\Omega)^*=\id_{\Omega^*}$, for all $\Phi,\Psi\in\Hom(\Mf_\Zs)$ and $\Omega\in\Ob(\Mf_\Zs)$. The evaluation transform $\ev:\Ob(\Mf_\Zs)\to\Hom(\Mf_\Zs)$, given by $\Omega\mapsto \ev^\Omega$ (where $\ev^\Omega:\Omega\to\Omega^{**}$ is the $\Zs$-linear map $x\mapsto \ev^\Omega_x$ that, for $x\in\Omega$, is defined as $\ev_x^\Omega(\phi):=\phi(x)$, for all $\phi\in\Omega$) is a natural tranformation between the covariant functors $\id_{\Mf_\Zs}\xrightarrow{\ev}*\circ *$, since: for all morphisms $\Omega_1\xrightarrow{\Phi}\Omega_2$ in $\Mf_\Zs$, we have $\ev^{\Omega_2}\circ\ \Phi=\Phi^{**}\circ\ev^{\Omega_1}$.  

\medskip 

We observe that the natural transformation $\id_{\Mf_\Zs}\xrightarrow{\ev}*\circ*$ satisfies the following ``weakened version'' of the triangle (right-right contravariant) adjunction identities: 
\begin{itemize}
\item 
$(\ev^\Omega)^*\circ \ev^{\Omega^*}=\id_{\Omega^*}$, for all $\Omega\in\Ob(\Mf_\Zs)$,  
\item
$\ev^{\Omega^*}\circ\ (\ev^\Omega)^*:\Omega^{***}\to\Omega^{***}$, for all $\Omega\in\Ob(\Mf_\Zs)$, is an idempotent ``projecting'' the $\Zs$-bimodule $\Omega^{***}$ onto its $\Zs$-submodule 
$\ev^{\Omega^*}(\Omega^*):=\{\ev^{\Omega^*}_\phi \ | \ \phi\in\Omega^* \}\subset\Omega^{***}$. 
\end{itemize}
Whenever we impose on the objects $\Omega$, the condition of reflexivity, i.e.~we ask that  $\Omega\xrightarrow{\ev^\Omega}\Omega^{**}$ is an isomorphism in $\Mf_\Zs$, the previous contravariant right semi-adjunction becomes a duality (contravariant right-right adjoint equivalence). 
\end{proof}

In the case of $\Zs$-central multimodules over $\Zs$-central $\Rs$-algebras, proposition~\ref{prop: sadj-z} generalizes in the form described in theorem~\ref{th: sadj}. 

\section{First-Order Differential Operators on $\Zs$-central Multimodules}\label{sec: 1st-ord-multi}

In this last appendix, we briefly preview a definition of first-order differential operator between multimodules. 

The first-order condition~\ref{def: 1st-ord}, is just a reformulation for multimodules of the usual first-order condition for operators acting on $\As$-bimodules, put forward in~\cite[sections $4.\gamma$ and $4.\delta$]{Co94}. 
We encountered such notion during our investigation of non-commutative vector fields and derivations of non-commutative algebras (see footnote~\ref{foo: ncdc}) and although we needed there to consider mostly first-order differential operators defined on $\As\otimes_\Zs\As$, for a non-commutative unital associative $\Rs_\Zs$-algebra $\As$, here for completeness we present the basic definition in a more general context. 

\begin{definition}\label{def: 1st-ord}
Let ${}_{(\As_\alpha)}\Ms_{(\Bs_\beta)}$ and ${}_{(\As_\alpha)}\Ns_{(\Bs_\beta)}$ be two multimodules, over the families of unital associative \hbox{$\Rs_\Zs$-algebras} $(\As_\alpha)_{\alpha\in A}$ and $(\Bs_\beta)_{\beta\in B}$.  
The set $\Diff^1_{(\As_\alpha)-(\Bs_\beta)}(\Ms;\Ns)$ of \emph{first-order differential operators} from $\Ms$ to $\Ns$ consists of those $\Rs$-linear maps $\Ms\xrightarrow{\delta}\Ns$ that satisfy the first-order conditions, for all $x\in\Ms$:  
\begin{gather}
\label{eq: 1st-ord}
\delta(a \cdot_\alpha x \cdot_\beta b) + a \cdot_\alpha \delta(x) \cdot_\beta  b
= \delta(a \cdot_\alpha x) \cdot_\beta b + a\cdot_\alpha \delta(x\cdot_\beta b), 
\quad 
\forall (\alpha,\beta)\in A\times B, \ 
(a,b)\in\As_\alpha\times\Bs_\beta,    
\\ \label{eq: 1st-ord-b}
\delta(a \cdot_\alpha b \cdot_{\alpha'} x) + a \cdot_\alpha b \cdot_{\alpha'}  \delta(x) 
= b \cdot_{\alpha'} \delta(a \cdot_\alpha x)  + a\cdot_\alpha \delta(b\cdot_{\alpha'} x), 
\ 
\forall \alpha\neq\alpha'\in A, \  
(a,b)\in\As_\alpha\times\As_{\alpha'}, 
\\ \label{eq: 1st-ord-c}
\delta(x \cdot_\beta a \cdot_{\beta'} b) + \delta(x) \cdot_\beta a \cdot_{\beta'} b   
= \delta(x\cdot_\beta a)  \cdot_{\beta'} b + \delta(x\cdot_{\beta'} b)\cdot_\beta a, 
\quad 
\forall \beta\neq\beta'\in B, \ 
(a,b)\in\Bs_\beta\times\Bs_{\beta'}. 
\end{gather}
\end{definition}

\begin{remark}
Defining $L^\alpha_a(x):=a\cdot_\alpha x$ and $R^\beta_b(x):=x\cdot_\beta b$, for all $(\alpha,\beta)\in A\times B$, $(a,b)\in\As_\alpha\times\Bs_\beta$ and $x\in\Ms$, 
a direct computation assures that, for all $(\alpha,\beta)\in A\times B$, equation~\eqref{eq: 1st-ord} above is equivalent to
\begin{equation*}
[ [\delta,L^\alpha_a]_-,R^\beta_b]_-=0_{\Hom_\Zs(\Ms;\Ns)}=[ [\delta,R^\beta_b]_-,L^\alpha_a]_-, 
\quad \forall (a,b)\in \As_\alpha\times\Bs_\beta, 
\end{equation*}
that is the familiar Connes' first-order condition for the operator $\Ms\xrightarrow{\delta}\Ns$ on the bimodules ${}_{\As_\alpha}\Ms_{\Bs_\beta}$, ${}_{\As_\alpha}\Ns_{\Bs_\beta}$. 

\medskip 

A perfectly similar reinterpretation in terms of commutators, provides: 
\begin{gather*}
[[\delta,L^\alpha_a]_-,L^{\alpha'}_b]_-=0_{\Hom_\Zs(\Ms;\Ns)}=[ [\delta,L^{\alpha'}_b]_-,L^\alpha_a]_-, 
\quad \forall \alpha\neq\alpha'\in A, \ (a,b)\in\As_\alpha\times\As_{\alpha'}, 
\end{gather*}
as an equivalent reformulation of equation~\eqref{eq: 1st-ord-b} for the left-$(\As_\alpha,\As_{\alpha'})$ bimodules  ${}_{\As_\alpha,\As_{\alpha'}}\Ms$, ${}_{\As_\alpha,\As_{\alpha'}}\Ns$ and 
\begin{gather*}
[[\delta,R^\beta_a]_-,R^{\beta'}_b]_-=0_{\Hom_\Zs(\Ms;\Ns)}=[ [\delta,R^{\beta'}_b]_-,R^\beta_a]_-,
\quad \forall \beta\neq\beta'\in B, \ (a,b)\in\As_\beta\times\As_{\beta'}, 
\end{gather*}
as a replacement of equation~\eqref{eq: 1st-ord-c} for the right-$(\Bs_\beta,\Bs_{\beta'})$ bimodules $\Ms_{\Bs_\beta,\Bs_{\beta'}}$, $\Ns_{\Bs_\beta,\Bs_{\beta'}}$. 
\xqed{\lrcorner}
\end{remark}

\begin{remark}
Definition~\ref{def: 1st-ord} is actually a special ($\Rs$-linear covariant) case of a much more general notion of first-order differential operator that allows to discuss $\Zs$-linear differential operators that are possibly contravariant and $\Rs$-conjugate linear. Making use of exactly the same notations introduced in definition~\ref{def: morphism} for morphisms (zero-order differential operators) between multimodules, we say that a \emph{first-order differential operator between multimodules} consists of the following data and conditions: 
\begin{gather*}
{}_{(\As_\alpha)_A}\Ms_{(\Bs_\beta)_B}\xrightarrow{(\phi,\eta,\delta,\zeta,\psi)_f}{}_{(\Cs_\gamma)_C}\Ns_{(\Ds_\delta)_D},
\quad 
A\uplus B\xrightarrow{f} C\uplus D,\quad A\xrightarrow{\phi}\End_\Zs(\Rs)\xleftarrow{\psi}B, \quad  
\delta\in\Hom_\Zs(\Ms;\Ns), 
\\
\forall (\alpha,\beta)\in A_+\times B_+, 
\quad \As_\alpha\xrightarrow{(\phi_\alpha,\eta_\alpha)}\Cs_{f(\alpha)},\  \Bs_\beta\xrightarrow{(\psi_\beta,\zeta_\beta)}\Ds_{f(\beta)},  
\ \text{$\Zs$-linear covariant unital homomorphisms},
\\
\forall (\alpha,\beta)\in A_-\times B_-, \quad 
\As_\alpha\xrightarrow{(\phi_\alpha,\eta_\alpha)}\Ds_{f(\alpha)}, \ \Bs_\beta\xrightarrow{(\psi_\beta,\zeta_\beta)}\Cs_{f(\beta)}, \ 
\text{$\Zs$-linear contravariant unital homomorphisms},  
\\
\delta(r\cdot_\alpha x\cdot_\beta s)=\phi_\alpha(r)\cdot_{f(\alpha)} \delta(x)\cdot_{f(\beta)} \psi_\beta(s), 
\quad \forall (\alpha,\beta)\in A_+\times B_+, \ r,s\in\Rs, \ x\in\Ms,
\\ 
\delta(r\cdot_\alpha x\cdot_\beta s)=\psi_\beta(s)\cdot_{f(\beta)} \delta(x)\cdot_{f(\alpha)} \phi_\alpha(r), 
\quad \forall (\alpha,\beta)\in A_-\times B_-, \ r,s\in\Rs, \ x\in\Ms,
\\
\forall (\alpha,\beta)\in A_+\times B_+, \  
\delta(a \cdot_\alpha x \cdot_\beta b) + \eta_\alpha(a) \cdot_{f(\alpha)} \delta(x) \cdot_{f(\beta)}  \zeta_\beta(b)
= \delta(a \cdot_\alpha x) \cdot_{f(\beta)} \zeta(b)_\beta + \eta_\alpha(a)\cdot_{f(\alpha)} \delta(x\cdot_\beta b), 
\\ 
\forall (\alpha,\beta)\in A_-\times B_-, \  
\delta(a \cdot_\alpha x \cdot_\beta b) + \zeta_\beta(b)\cdot_{f(\beta)}\delta(x)\cdot_{f(\alpha)} \eta_\alpha(a)    
= \zeta(b)_\beta\cdot_{f(\beta)}\delta(a \cdot_\alpha x)   + \delta(x\cdot_\beta b)\cdot_{f(\alpha)}\eta_\alpha(a), 
\\ 
\forall (\alpha,\beta)\in A_+\times B_-, \  
\delta(a \cdot_\alpha x \cdot_\beta b) +  \zeta_\beta(b)\cdot_{f(\beta)}\eta_\alpha(a) \cdot_{f(\alpha)} \delta(x)  
=  \zeta(b)_\beta\cdot_{f(\beta)}\delta(a \cdot_\alpha x)  + \eta_\alpha(a)\cdot_{f(\alpha)} \delta(x\cdot_\beta b), 
\\ 
\forall (\alpha,\beta)\in A_-\times B_+, \  
\delta(a \cdot_\alpha x \cdot_\beta b) +  \delta(x) \cdot_{f(\beta)}  \zeta_\beta(b)\cdot_{f(\alpha)}\eta_\alpha(a) 
= \delta(a \cdot_\alpha x) \cdot_{f(\beta)} \zeta(b)_\beta +  \delta(x\cdot_\beta b)\cdot_{f(\alpha)}\eta_\alpha(a), 
\\ 
\forall \alpha,\alpha'\in A_+, \ 
\delta(a \cdot_\alpha  b\cdot_{\alpha'} x) + \eta_\alpha(a) \cdot_{f(\alpha)} \eta_{\alpha'}(b) \cdot_{f(\alpha')}\delta(x) 
= \eta_{\alpha'}(b) \cdot_{f(\alpha')}\delta(a \cdot_\alpha x) + \eta_\alpha(a)\cdot_{f(\alpha)} \delta(b\cdot_{\alpha'} x), 
\\ 
\forall \alpha,\alpha'\in A_-, \ 
\delta(a \cdot_\alpha  b\cdot_{\alpha'} x) + \delta(x) \cdot_{f(\alpha)}\eta_\alpha(a)  \cdot_{f(\alpha')}\eta_{\alpha'}(b) 
=  \delta(a \cdot_\alpha x)\cdot_{f(\alpha')}\eta_{\alpha'}(b) + \delta(b\cdot_{\alpha'} x)\cdot_{f(\alpha)} \eta_\alpha(a), 
\\ 
\forall \beta,\beta'\in B_+, \ 
\delta(x \cdot_\beta a  \cdot_{\beta'}b ) + \delta(x)\cdot_{f(\beta)} \zeta_\beta(a) \cdot_{f(\beta')}\zeta_{\beta'}(b)  
=  \delta(x\cdot_\beta a )\cdot_{f(\beta')}\zeta_{\beta'}(b) + \delta(x\cdot_{\beta'} b)\cdot_{f(\beta)} \zeta_\beta(a), 
\\ 
\forall \beta,\beta'\in B_-, \ 
\delta(x \cdot_\beta a  \cdot_{\beta'}b ) +  \zeta_\beta(a) \cdot_{f(\beta)}\zeta_{\beta'}(b) \cdot_{f(\beta')}\delta(x) 
= \zeta_{\beta'}(b) \cdot_{f(\beta')}\delta(x\cdot_\beta a ) + \zeta_\beta(a)\cdot_{f(\beta)} \delta(x\cdot_{\beta'} b), 
\\ 
\forall \alpha_\pm\in A_\pm,\ 
\delta(a \cdot_{\alpha_+}  b\cdot_{\alpha_-} x) + \eta_{\alpha_+}(a) \cdot_{f(\alpha_+)}  \delta(x) \cdot_{f(\alpha_-)}\eta_{\alpha_-}(b)
=  \delta(a \cdot_\alpha x)\cdot_{f(\alpha_-)}\eta_{\alpha_-}(b) + \eta_\alpha(a)\cdot_{f(\alpha)} \delta(b\cdot_{\alpha_-} x), 
\\ 
\forall \beta_\pm\in B_\pm, \  
\delta(x \cdot_{\beta_+} a \cdot_{\beta_-}b) +   \zeta_{\beta_-}(b)\cdot_{f(\beta_-)}\delta(x)\cdot_{f(\beta_+)}\zeta_{\beta_+}(a) 
=  \zeta_{\beta_-}(b)\cdot_{f(\beta_-)}\delta(x \cdot_{\beta_+} a) +  
\delta(x\cdot_{\beta_-} b)\cdot_{f(\beta_+)}\zeta_{\beta_+}(a).  
\end{gather*}
Whenever $f$, $(\phi_\alpha)_A$, $(\psi_\beta)_B$, $(\eta_\alpha)_A$ and $(\zeta_\beta)_B$ are all identity functions, we recover the initial definition~\ref{def: 1st-ord}. 

\medskip 

Given an arbitrary signature $\sigma:=(\phi,\eta,\zeta,\psi)_f$, we denote by $\Diff^1_\sigma(\Ms;\Ns)$ the family of first-order differential operators defined by all the conditions above. Making use of remark~\ref{prop: cj-duals} we can obtain a bijective correspondence between $\Diff^1_\sigma(\Ms;\Ns)$ and $\Diff^1_{(\As_\alpha)_A-(\Bs_\beta)_B}(\Ms;\Ns^\sigma)$ that is associating to each first-order differential operator $\delta\in\Diff^1_\sigma(\Ms;\Ns)$, with signature $\sigma$, the unique first-order differential operator $\delta^\sigma\in\Diff^1_{(\As_\alpha)_A-(\Bs_\beta)_B}(\Ms;\Ns^\sigma)$ such that $\Theta^\sigma_\Ns\circ\delta^\sigma=\delta$.
\xqed{\lrcorner}
\end{remark}

The family of first-order differential operators $\Diff^1_{(\As_\alpha)_A-(\Bs_\beta)_B}(\Ms;\Ns)$ is a central $\Zs$-bimodule, but does not usually have other well-defined actions, even of the algebra $\Rs$. 
Whenever all the $\Rs_\Zs$-algebras involved are $\Rs$-central bimodules (in particular if $\Rs=\Zs$) the following immediate result is of interest. 

\begin{remark}
Let ${}_{(\As_\alpha)_A}\Ms_{(\Bs_\beta)_B}$ and ${}_{(\As_\alpha)_A}\Ns_{(\Bs_\beta)_B}$ be multimodules over $\Rs_\Zs$-central algebras. 

For any pair of sub-families of indexes $I\times J\subset A\times B$, define 
$\Is^I_\alpha:=
\begin{cases} 
\As_\alpha, \quad \alpha\in I 
\\ 
\Rs,\phantom{_\alpha} \quad \alpha\notin I
\end{cases},
\quad \Js^J_\beta:=
\begin{cases} 
\Bs_\beta, \quad \beta\in J 
\\ 
\Rs,\phantom{_\beta} \quad \beta\notin J
\end{cases}$.

The spaces $\Diff^1_{(\Is^I_\alpha)-(\Js^J_\beta)}(\Ms;\Ns)$ are all $\Zs$-central multimodules with respect to the following actions: 
\begin{gather*}
(a\cdot_\alpha \delta \cdot_\beta b) (x):=a\cdot_\alpha\delta(x)\cdot_\beta b, \quad \forall (\alpha,\beta)\in (A-I)\times (B-J), \ \ (a,b)\in \Is^I_\alpha\times\Js^J_\beta=\As_\alpha\times\Bs_\beta, \  x\in\Ms, 
\\
(b\odot_\beta \delta \odot_\alpha a) (x):=\delta(a\cdot_\alpha x\cdot_\beta b), \quad \forall (\alpha,\beta)\in (A-I)\times (B-J), \ \ (a,b)\in \Is^I_\alpha\times\Js^J_\beta=\As_\alpha\times\Bs_\beta, \  x\in\Ms, 
\\
(r\cdot_\alpha \delta \cdot_\beta s) (x):=r\cdot_\alpha\delta(x)\cdot_\beta s, \quad \forall (\alpha,\beta)\in I\times J, \ \ (r,s)\in \Is^I_\alpha\times\Js^J_\beta=\Rs, \  x\in\Ms, 
\\
(s\odot_\beta \delta \odot_\alpha r) (x):=\delta(r\cdot_\alpha x\cdot_\beta s), \quad \forall (\alpha,\beta)\in I\times J, \ \ (r,s)\in \Is^I_\alpha\times\Js^J_\beta=\Rs, \  x\in\Ms.
\end{gather*}
Whenever $(I_1,J_1)\leq (I_2,J_2)$, we have inclusions $\Diff^1_{(\Is^{I_2}_\alpha)_A-(\Js^{J_2}_\beta)_B}(\Ms;\Ns) \subset\Diff^1_{(\Is^{I_1}_\alpha)_A-(\Js^{J_1}_\beta)_B}(\Ms;\Ns)$ of $\Zs$-central bimodules that are also a morphisms of $\Zs$-central multimodules with respect to all the relevant actions. 
\xqed{\lrcorner}
\end{remark}

\end{document}